\documentclass[vecarrow]{svmult}

\usepackage{makeidx}         
\usepackage{graphicx}        
\usepackage{psfrag}          
\usepackage{multicol}        
\usepackage{amsmath}
\usepackage{amssymb}
\usepackage{amsfonts}
\usepackage{latexsym}
\usepackage[T1]{fontenc}
\usepackage{subfigure}                      
\usepackage{float}
\usepackage{mathrsfs}
\usepackage{mathptmx}         

\makeindex   

\begin{document}

\title{Shape optimization for interior Neumann and transmission eigenvalues}
\index{shape optimization}
\index{interior Neumann eigenvalues}
\index{interior transmission eigenvalues}

\titlerunning{Shape optimization}

\author{A. Kleefeld}

\authorrunning{A. Kleefeld}

\institute{A. Kleefeld\at Forschungszentrum J\"ulich GmbH, Supercomputing Centre J\"ulich, 52425 J\"ulich, Germany,\hfill\break
\email{a.kleefeld@fz-jueliche.de}
}

\maketitle

\index{integral equations}
\index{non-linear eigenvalue problem}

\bgroup


\allowdisplaybreaks

\abstract{
Shape optimization problems for interior eigenvalues is a very challenging task since already the 
computation of interior eigenvalues for a given shape is far 
from trivial. 
For example, a concrete maximizer with respect to shapes of fixed area is theoretically established only for the first two non-trivial 
Neumann eigenvalues. The existence of such a maximizer for higher Neumann 
eigenvalues is still unknown. Hence, the problem should be addressed 
numerically. Better numerical results are achieved for the maximization of 
some Neumann eigenvalues using boundary integral equations for a simplified 
parametrization of the boundary in combination with a non-linear eigenvalue 
solver. 
Shape optimization for interior transmission eigenvalues is even more complicated 
since the corresponding transmission problem is non-self-adjoint and non-elliptic.
For the first time numerical results are presented for the 
minimization of interior transmission eigenvalues for which no single theoretical result is 
yet available.
}

\section{Introduction}\label{sec:kle1}
The task is to optimize the shape of a domain $\Omega\subset \mathbb{R}^2$ with respect to the $k$-th eigenvalue under the constraint that the area $|\Omega|$ of the 
domain is constant, say $A$. Here, the domain is an open and bounded set with smooth boundary $\partial \Omega$ which is also allowed to be disconnected.
In the sequel, we consider two different problems. 

First, we deal with the maximization of interior Neumann eigenvalues (INEs). Precisely, one has to
find numbers $\lambda>0$ such that 
\begin{eqnarray*}
\Delta u+\lambda u=0 \text{ in }\Omega\,,\qquad \partial_\nu u=0 \text{ on }\partial \Omega
\end{eqnarray*}
is satisfied for non-trivial $u$, where $\nu$ denotes the normal pointing in the exterior. It is well-known that this problem is elliptic and the eigenvalues are discrete. 
The case $\lambda=0$ which corresponds to a constant function is not considered here.
It has been shown in 1954 and 1956 that the first INE is maximized by a circle (see \cite{Sz54, We56}) and recently that the second INE is maximized by two disjoint circles of the same size (see \cite{GiNaPo09}). 
However, the existence and uniqueness of a shape maximizer for higher INEs is from the theoretically point of view still unknown.  
But numerical results suggest that such a maximizer might exist. We refer the reader to \cite{AnFr12,AnOu17} for recent results and a good overview over 
who has already worked in this direction. In Figure \ref{fig:kle1} we show numerically the shape maximizer for the first six INEs.
\begin{figure}
\centering
\begin{subfigure}{}
\includegraphics[width=3.5cm]{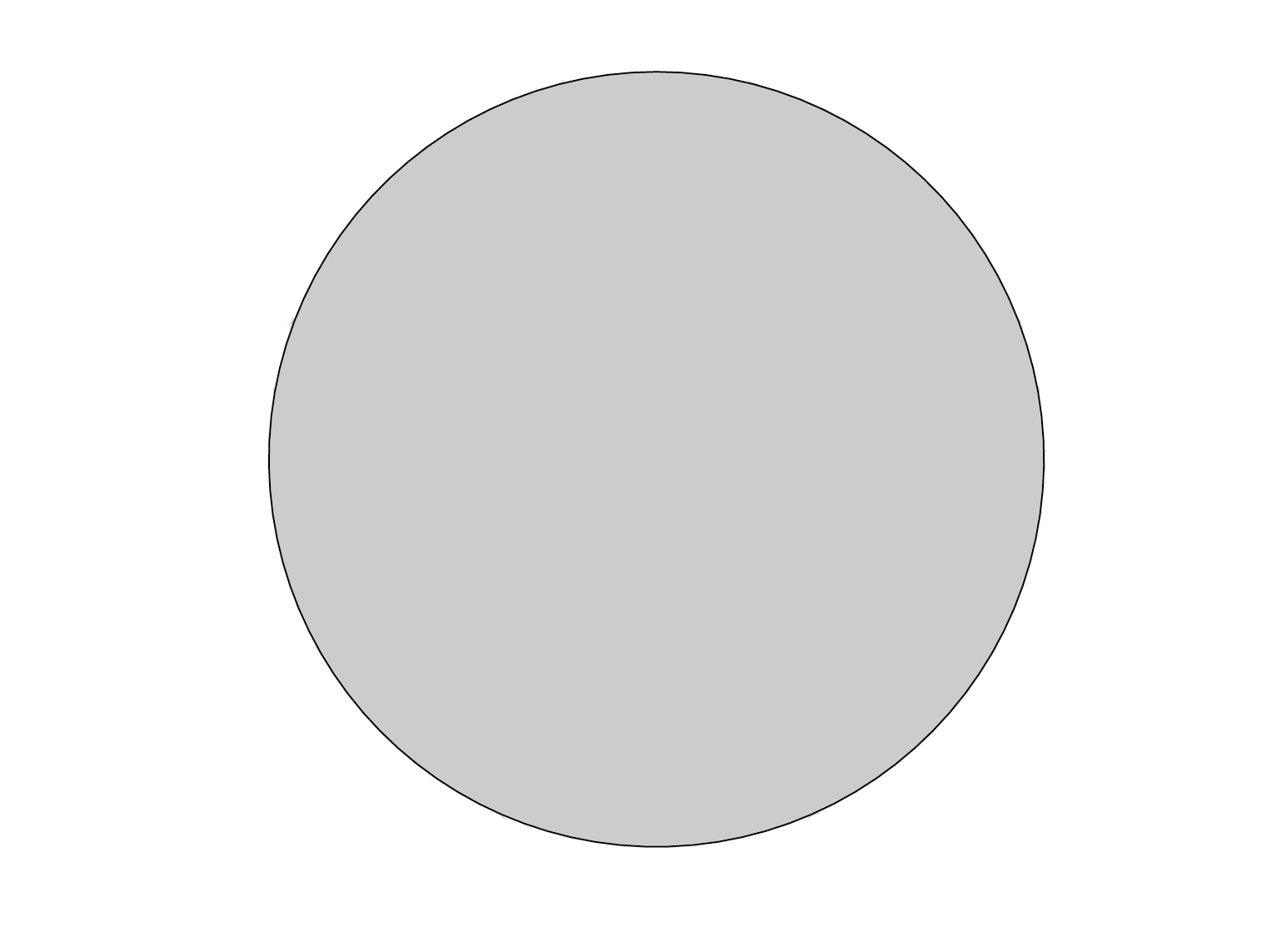}
\end{subfigure}\quad
\begin{subfigure}{}
\includegraphics[width=3.5cm]{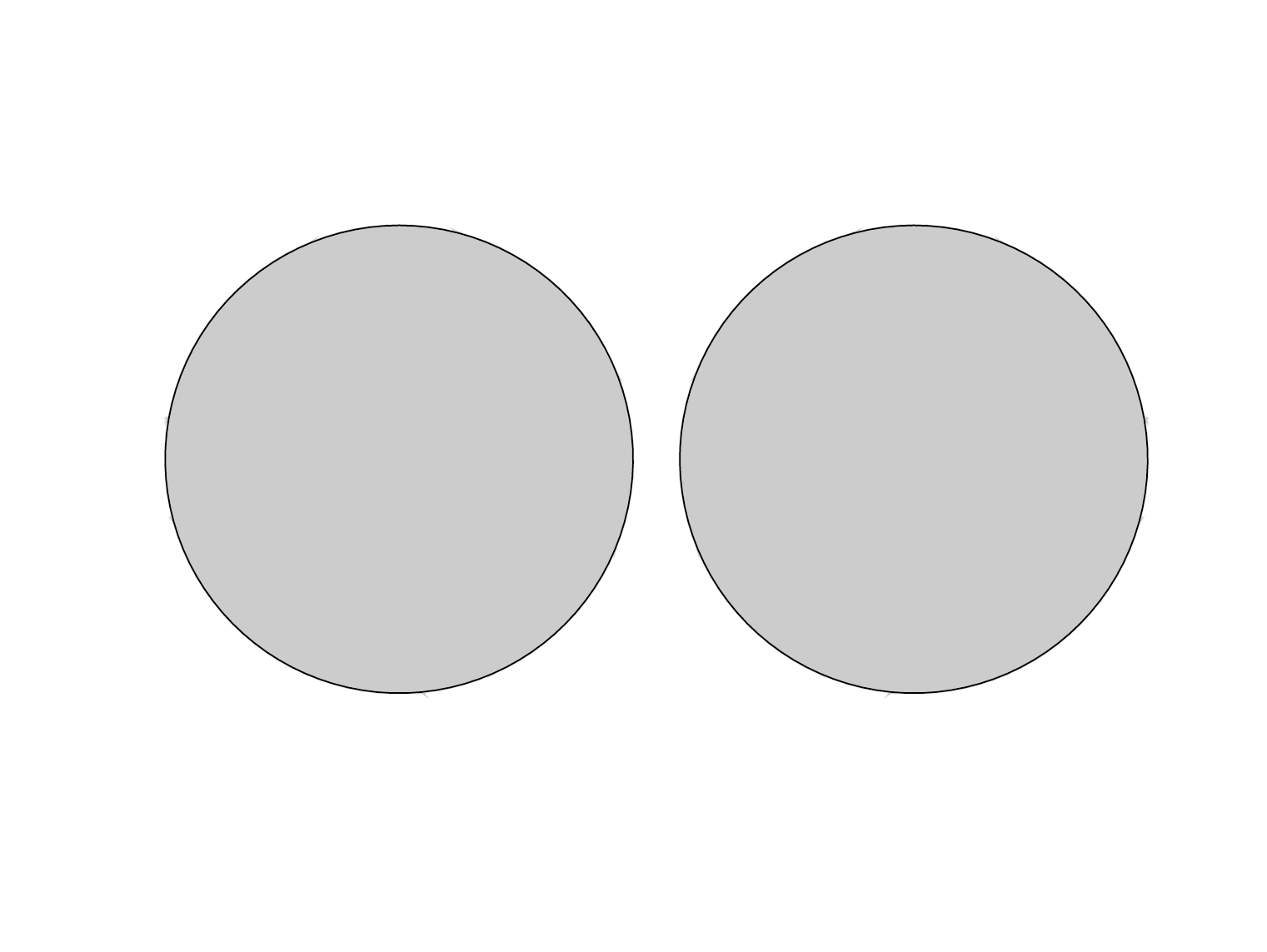}
\end{subfigure}\quad
\begin{subfigure}{}
\includegraphics[width=3.5cm]{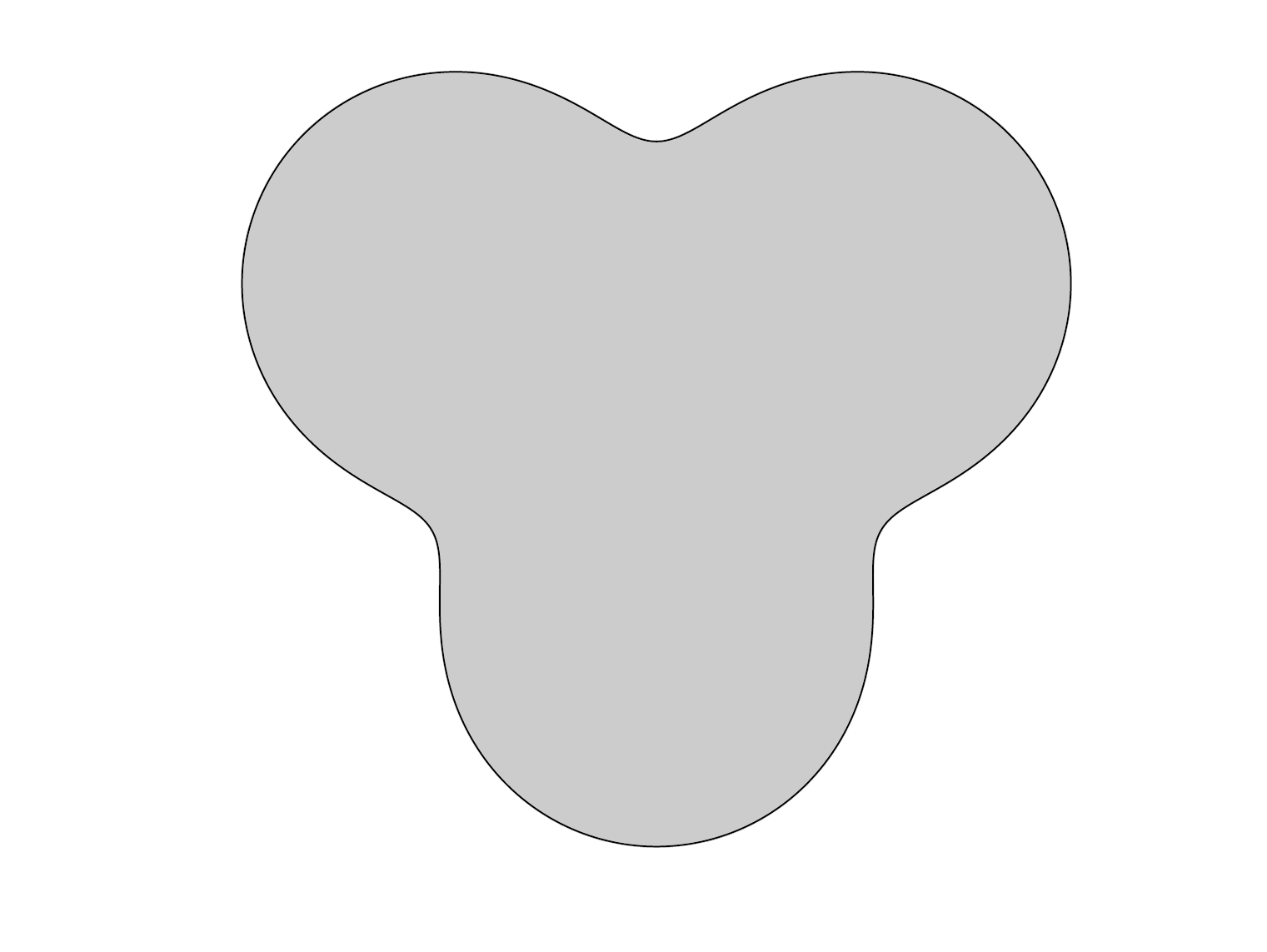}
\end{subfigure}\quad
\begin{subfigure}{}
\includegraphics[width=3.5cm]{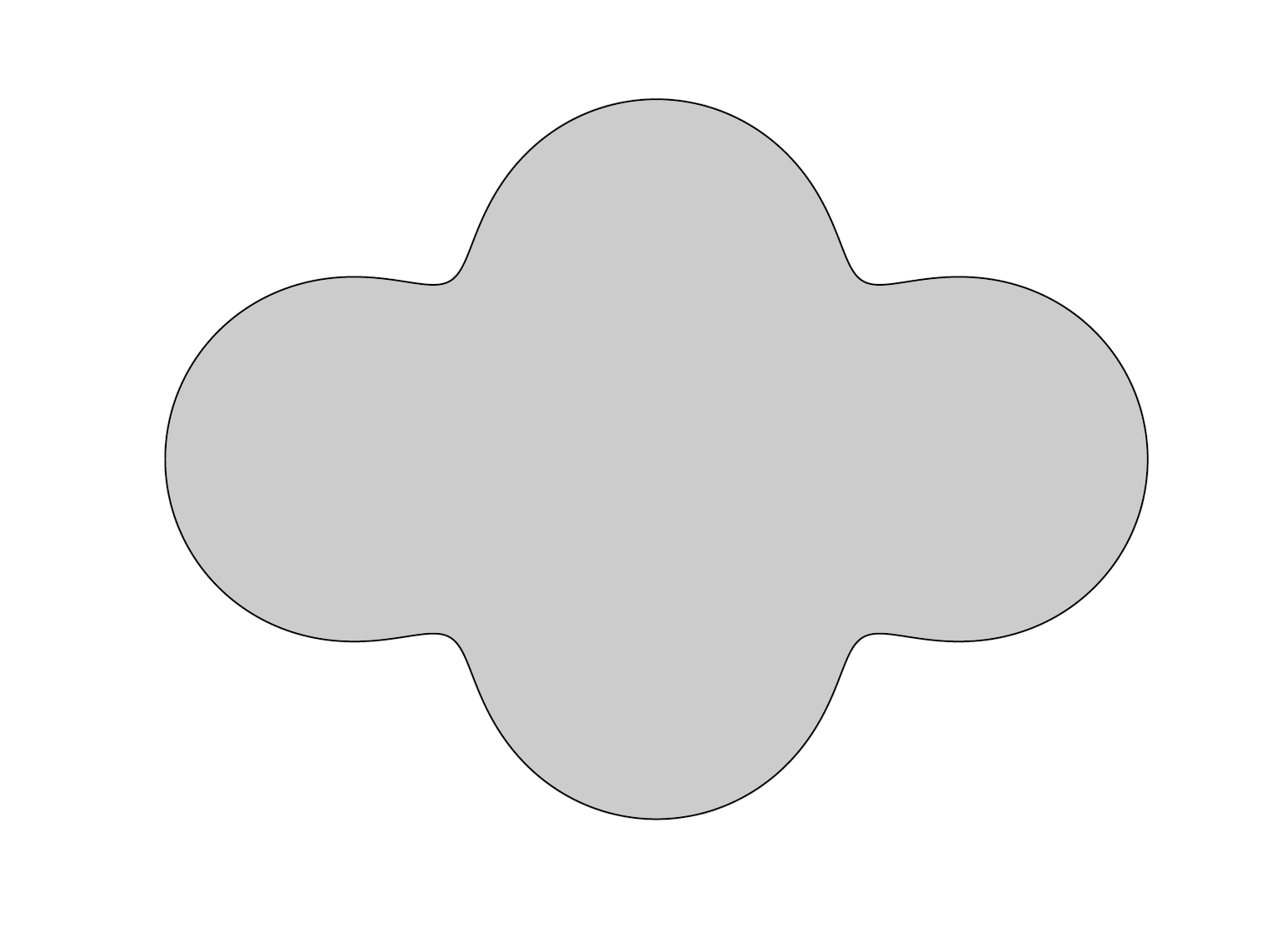}
\end{subfigure}\quad
\begin{subfigure}{}
\includegraphics[width=3.5cm]{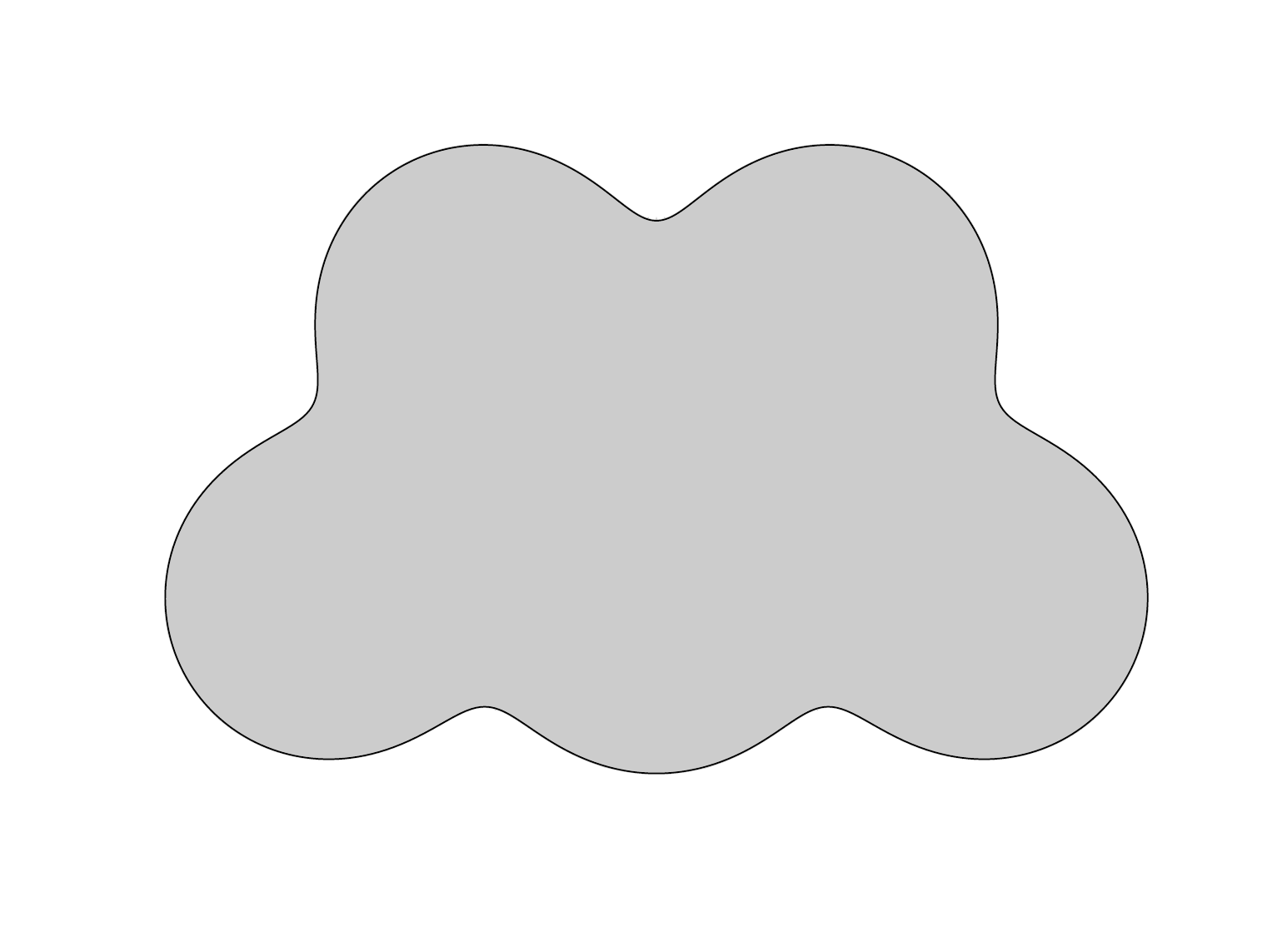}
\end{subfigure}\quad
\begin{subfigure}{}
\includegraphics[width=3.5cm]{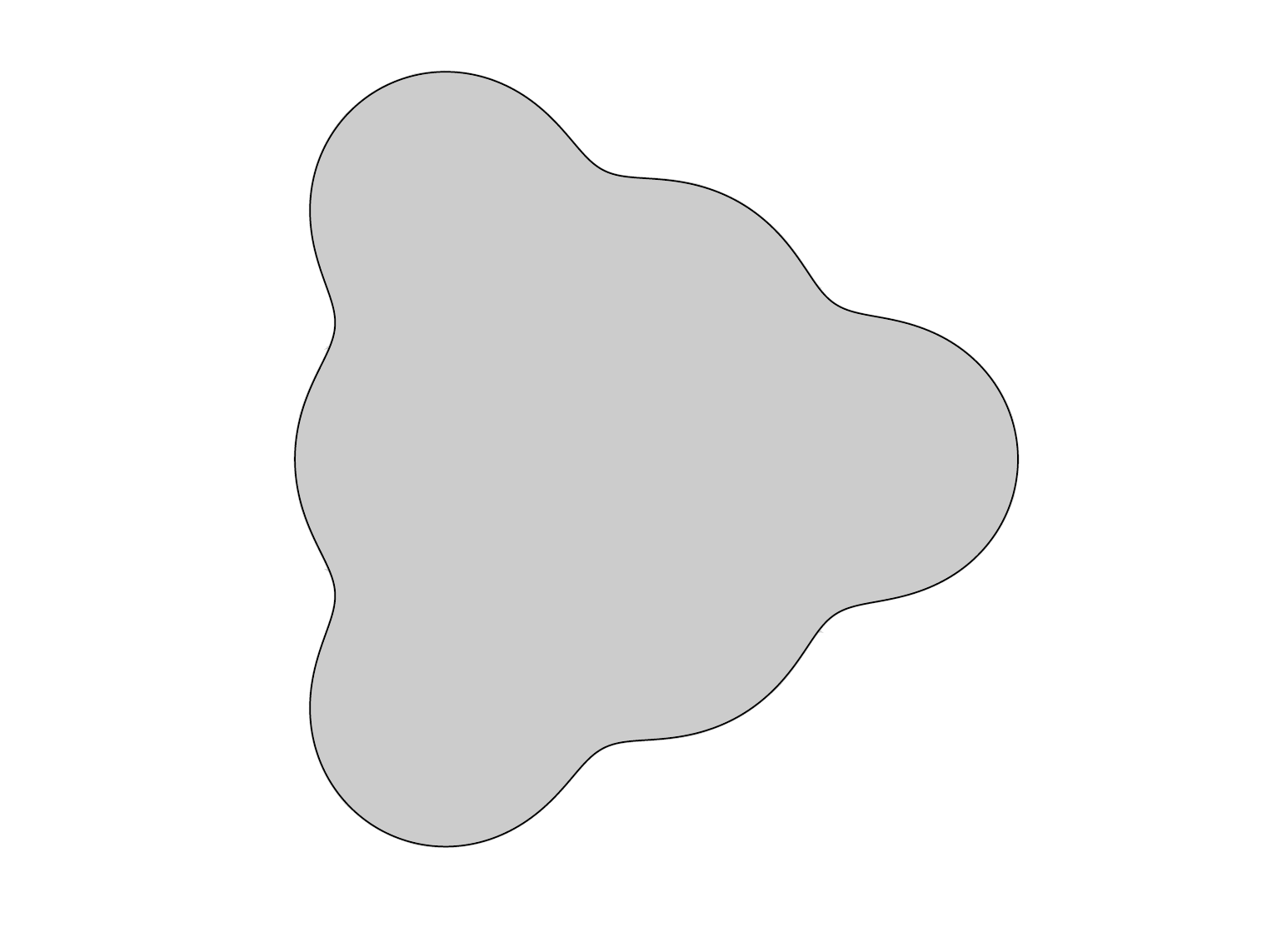}
\end{subfigure}
\caption{Shape maximizer for the first six INEs obtained numerically. The recent optimal values $\lambda_k\cdotp A$ for $k=1,\ldots,6$ are $10.66$, $21.28$, $32.90$, $43.86$, $55.17$, $67.33$ (see \cite{AnOu17}).}
\label{fig:kle1}
\end{figure}

The optimal values $\lambda_k\cdotp A$ for $k=1,\ldots,6$ are $10.66$, $21.28$, $32.79$, $43.43$, $54.08$, $67.04$ (see \cite{AnFr12}) which have been improved recently to 
$10.66$, $21.28$, $32.90$, $43.86$, $55.17$, $67.33$ (see \cite{AnOu17}). This paper reports improved values for the third and fourth INE and at the same time the boundary of the shape maximizer is described explicitly in terms of two parameters.   

The second problem under consideration is the interior transmission problem. Interior transmission eigenvalues (ITEs) 
are numbers $\lambda\in \mathbb{C}\backslash \{0\}$ such that 
    \begin{align*}
    \Delta w+\lambda n w&=0\quad\;\;\text{in}\; \Omega\;,\\
    \Delta v+\lambda\;\;\; v&=0\quad\;\;\text{in}\;\Omega\;,\\
                           v&=w\quad\;\text{on}\;\partial\Omega\;,\\
            \partial_{\nu}v&=\partial_{\nu}w\;\text{on}\;\partial\Omega\;,
    \end{align*}
has a non-trivial solution $(v,w)\neq(0,0)$, where $n$ is the given index of refraction. 
However, this is a non-elliptic and non-self-adjoint problem appearing first in 1986 (see \cite{Ki86}). 
Existence and discreteness for real-valued $\lambda$ has been shown in \cite{CaGiHa10}. But, the existence is still 
open for complex-valued $\lambda$ except for special geometries (see \cite{SlSt16,CoLe17}). The computation of ITEs for a given shape 
is therefore a very challenging task (see \cite{KlPi18} for an excellent overview of existing methods). It is also noteworthy that neither theoretical nor numerical results are available for a shape optimizer
of the first two ITEs. Within this paper we give numerical evidence for a shape minimizer of the first two ITEs and stating a conjecture which researcher in this field might want to prove in the future.

\subsection*{Contribution of the paper}
The contribution of this paper is twofold. First, improved numerical results for the maximization of some interior Neumann eigenvalues are presented using a simplified parametrization of the boundary. 
Second, the previous concept is transferred in order to obtain numerical results for the minimization of interior transmission eigenvalues for the first time for which no single theoretical result is 
yet available.

\subsection*{Outline of the paper}
The paper is organized as follows: In Section \ref{sec:kle2}, it is explained in detail how to compute interior Neumann eigenvalues using a boundary integral equation followed by its discretization. 
Then, it is described how the resulting non-linear eigenvalue problem is solved numerically. Further, the new parametrization is introduced and used to obtain improved numerical results for the maximization of some 
interior Neumann eigenvalues. In Section \ref{sec:kle3}, the concept of the previous section is applied for the minimization of interior transmission eigenvalues for which 
neither numerical results nor theoretical results are yet available. Finally, a short summary and an outlook is given in Section \ref{sec:kle4}. 

\section{Shape optimization for interior Neumann eigenvalues}\label{sec:kle2}
Recall that interior Neumann eigenvalues (INEs) are numbers $\lambda=\kappa^2$ such that 
\begin{eqnarray*}
\Delta u+\kappa^2 u=0 \text{ in }\Omega\,,\qquad \partial_\nu u=0 \text{ on }\partial \Omega
\end{eqnarray*}
is satisfied. Note that this problem is elliptic and it is well-known that the eigenvalues are discrete and positive real-valued numbers. In the sequel, we ignore $\kappa=0$ which corresponds to the constant function.
To find such INEs for a given domain $\Omega$, we use a boundary integral equation approach. A single layer ansatz with unknown density $\psi$ given by
\begin{eqnarray*}
u(X)=\int_{\partial \Omega} \Phi_{\kappa}(X,y) \psi(y)\,\mathrm{d}s(y)\,,\quad X\in\Omega
\end{eqnarray*}
is used, where $\Phi_{\kappa}(X,y)=\mathrm{i}\,H_0^{(1)}(\kappa\|X-y\|)/4$ is the fundamental solution of the Helmholtz equation. Taking the normal derivative, $\Omega \ni X\rightarrow x\in\partial \Omega$, and 
using the jump condition yields the following boundary integral equation of the second kind
\begin{eqnarray}
\frac{1}{2}\psi(x)+\underbrace{\int_{\partial \Omega}\partial_{\nu(x)} \Phi_{\kappa}(x,y)\psi(y)\,\mathrm{d}s(y)}_{K(\kappa)}=0\,.
\label{bon:kle}
\end{eqnarray}
Note that the operator $K(\kappa):H^{-1/2}(\partial \Omega)\rightarrow H^{-1/2}(\partial \Omega)$ is compact assuming a regular boundary (see \cite{Mc00}). Hence, $Z(\kappa)=I/2+K(\kappa)$ is Fredholm of index zero for 
$\kappa\in\mathbb{C}\backslash \mathbb{R}_{\le 0}$ and thus the theory of eigenvalue problems for holomorphic Fredholm operator-valued functions applies to $Z(\kappa)$.

The integral equation (\ref{bon:kle}) is discretized via the boundary element collocation method. Precisely, we subdivide the boundary into $n/2$ pieces, approximate it by quadratic interpolation (the approximated boundary is denoted by $\widetilde{\partial \Omega}$), 
and define on each piece a quadratic interpolation for $\psi$. This leads to 
\begin{eqnarray*}
\underbrace{\left(\frac{1}{2}\mathbf{I}+\mathbf{M}(\kappa)\right)}_{\mathbf{Z}(\kappa)\in \mathbb{C}^{n\times n}}\vec{\psi}=\vec{0}\,,
\end{eqnarray*}
where the matrix entries of $\mathbf{M}$ are numerically calculated with the Gauss-Kronrad quadrature (see \cite{KlLi12} for details in the three-dimenensional case).
The resulting non-linear eigenvalue problem of the form
\begin{eqnarray*}
\mathbf{Z}(\kappa)\vec{\psi}=\vec{0}
\end{eqnarray*}
is solved with the method of Beyn \cite{Be12}. This method can find all eigenvalues $\kappa$ including their multiplicities within any contour
$\mathcal{C}\subset \mathbb{C}$ which is based on Keldysh's theorem. Precisely, one integrates the resolvent over the given contour whereas the integral is 
approximated with the trapezoidal rule (see \cite{Be12} for more details).
Hence, we are now able to compute highly accurate INEs for a given shape $\Omega$. Next, it is explained how to choose a parametrization for the boundary of $\Omega$. 
The idea is to use an implicit curve rather than an explicit representation of the curve.
 Equipotentials are implicit curves of the form
\begin{eqnarray}
\sum_{i=1}^m \frac{1}{\|x-P_i\|}=c\,,
\label{easy:kle}
\end{eqnarray}
where the parameter $c$ and the centers $P_i$ are given. Here, $\|\cdotp\|$ denotes the Euclidean norm. Precisely, all points $x\in \mathbb{R}^2$ satisfying (\ref{easy:kle}) for given points $P_i$, $i=1,\ldots, m$ and 
parameter $c$ describe the implicit curve.
\begin{example}\label{ex1:kle}
We choose three points $(-\sqrt{3}/2,1/2)$, $(\sqrt{3}/2,1/2)$, $(0,-1)$ for $m=3$ and $(-3/2,0)$, $(3/2,0)$, $(0,-\sqrt{3}/2)$, $(0,\sqrt{3}/2)$ for $m=4$. The edge length of the following geometric shapes as shown in Figure \ref{fig:kle2} is $\sqrt{3}$.
\begin{figure}
\centering
\begin{subfigure}{}
\includegraphics[width=5cm]{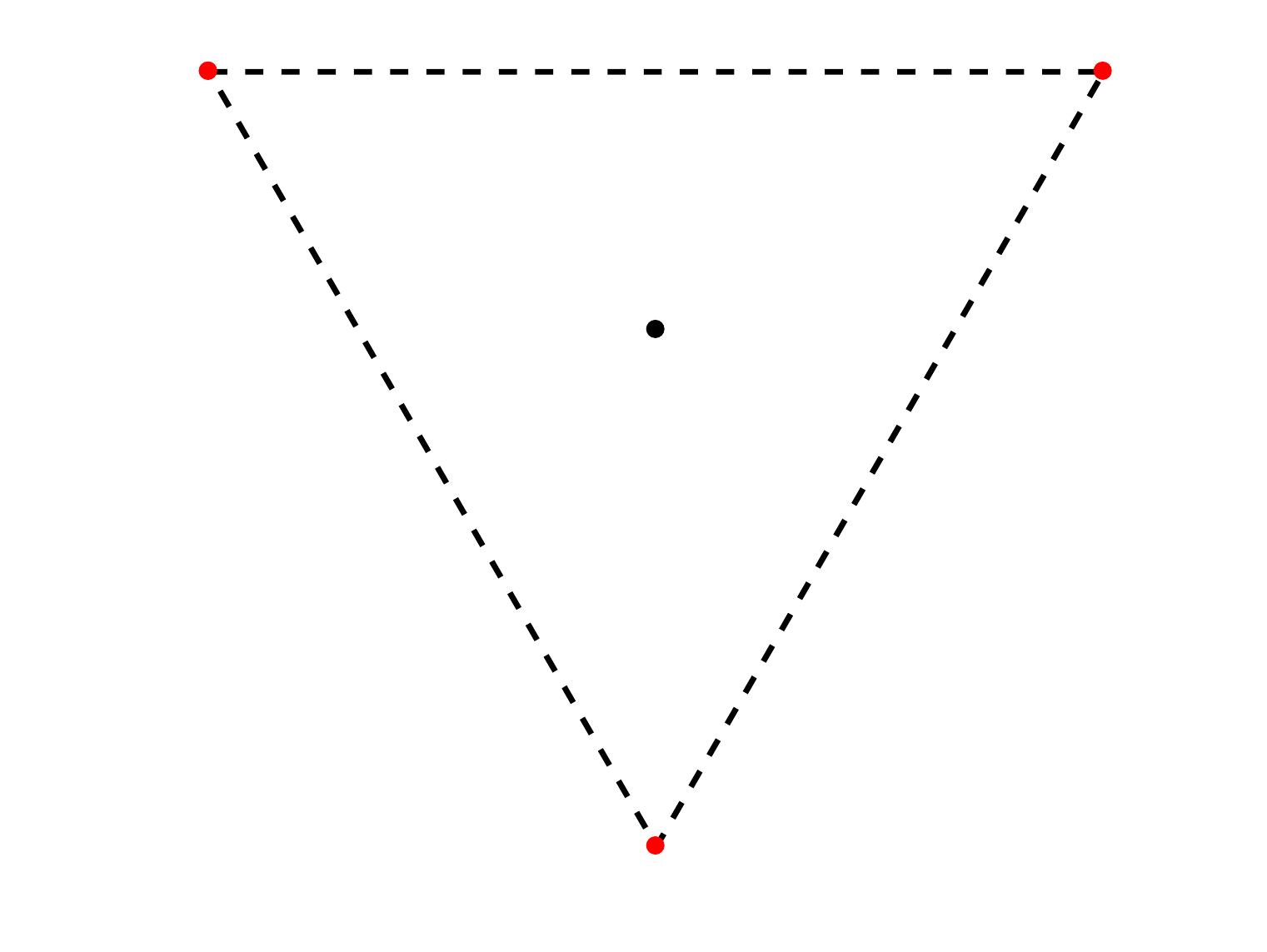}
\end{subfigure}\quad
\begin{subfigure}{}
\includegraphics[width=5cm]{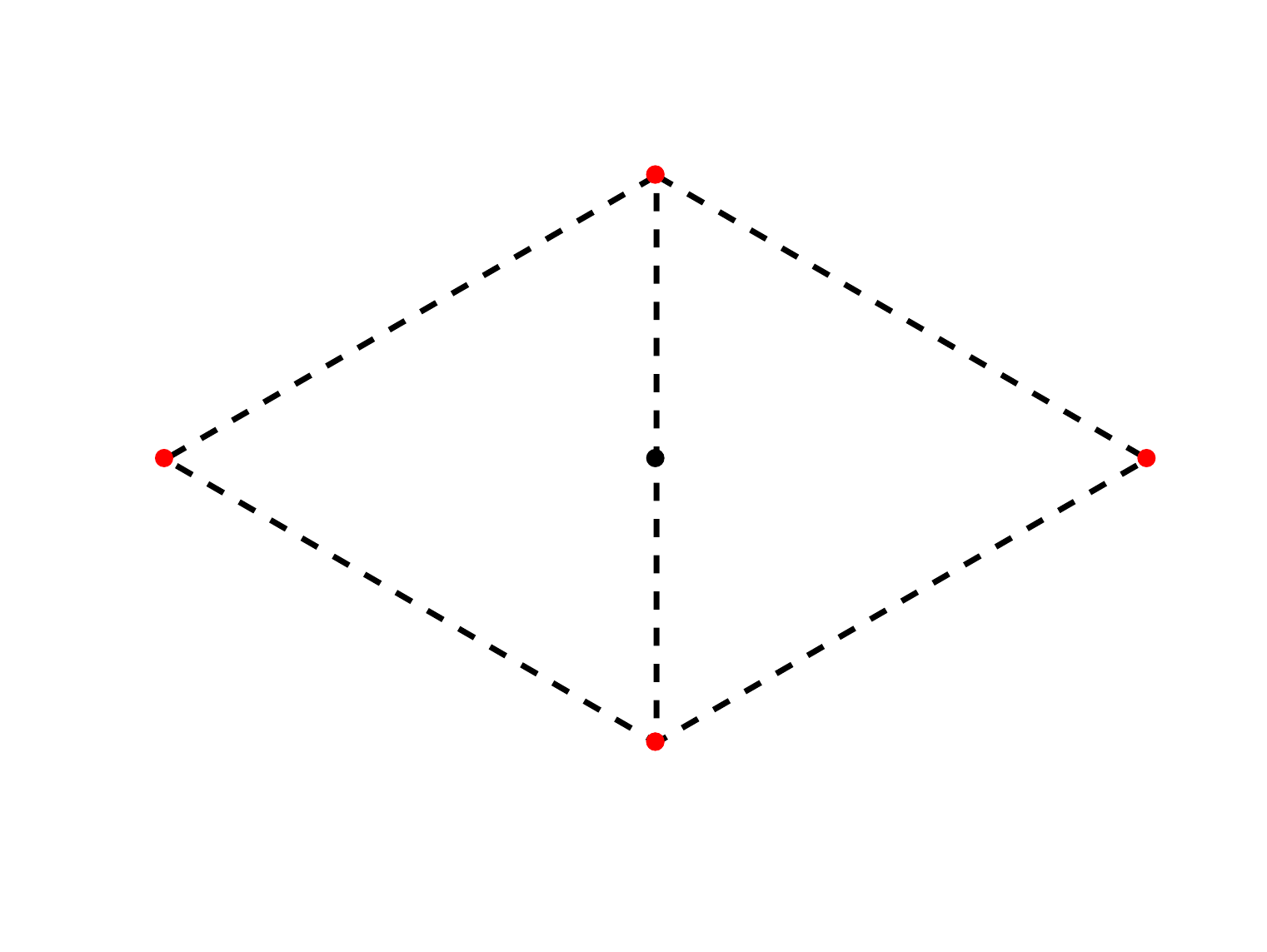}
\end{subfigure}
\caption{The choice of the points for $m=3$ are $(-\sqrt{3}/2,1/2)$, $(\sqrt{3}/2,1/2)$, $(0,-1)$ and for $m=4$ are $(-3/2,0)$, $(3/2,0)$, $(0,-\sqrt{3}/2)$, $(0,\sqrt{3}/2)$ shown as a red dot. 
The origin is shown as a black dot.}
\label{fig:kle2}
\end{figure}
\end{example}

Next, we show the influence of the parameter $c$. As one can see in Figure \ref{fig:kle3} the larger the parameter $c$ gets, the more constricting the boundary gets.
\begin{figure}
\centering
\begin{subfigure}{}
\includegraphics[width=3.5cm]{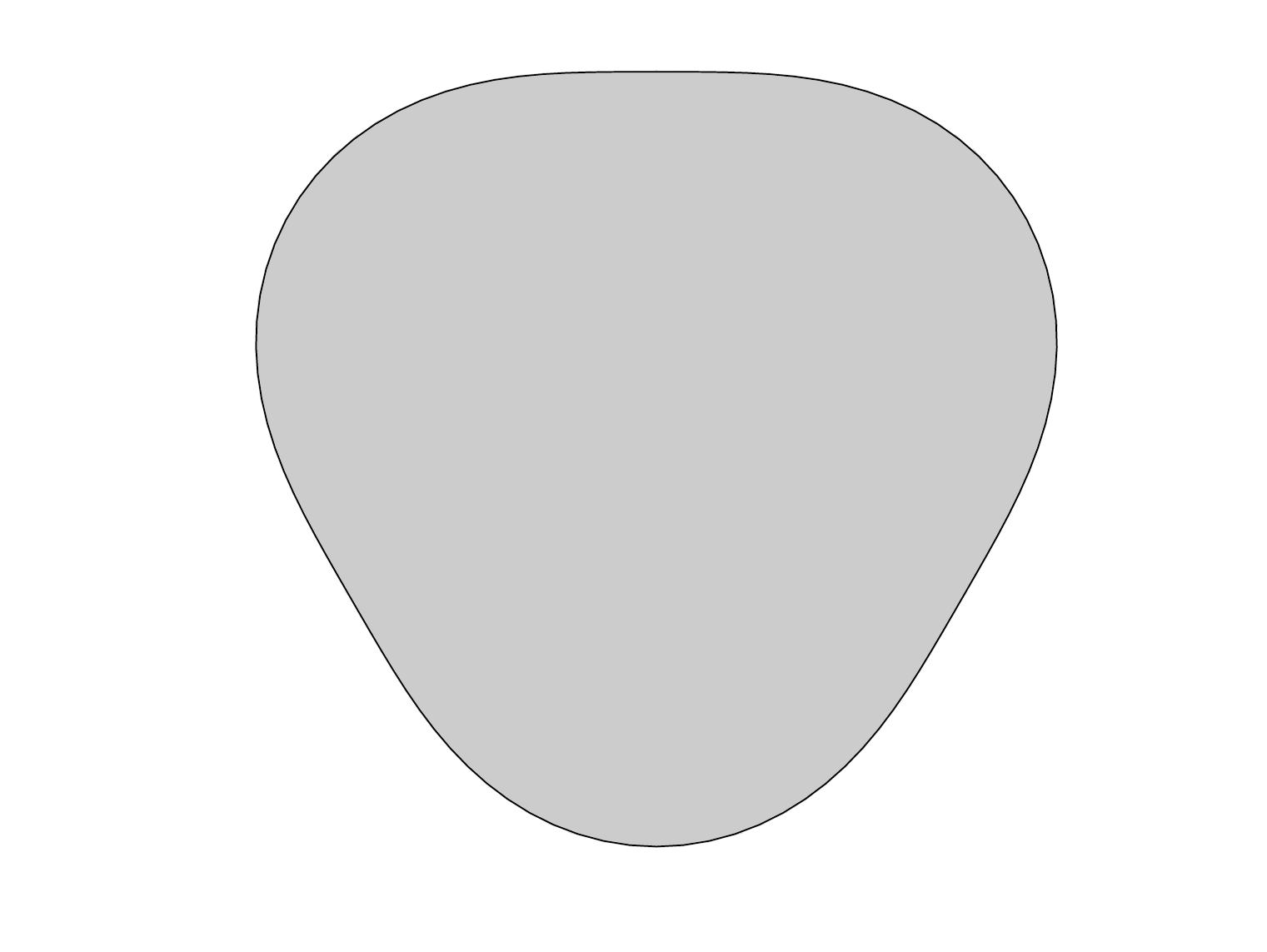}
\end{subfigure}\quad
\begin{subfigure}{}
\includegraphics[width=3.5cm]{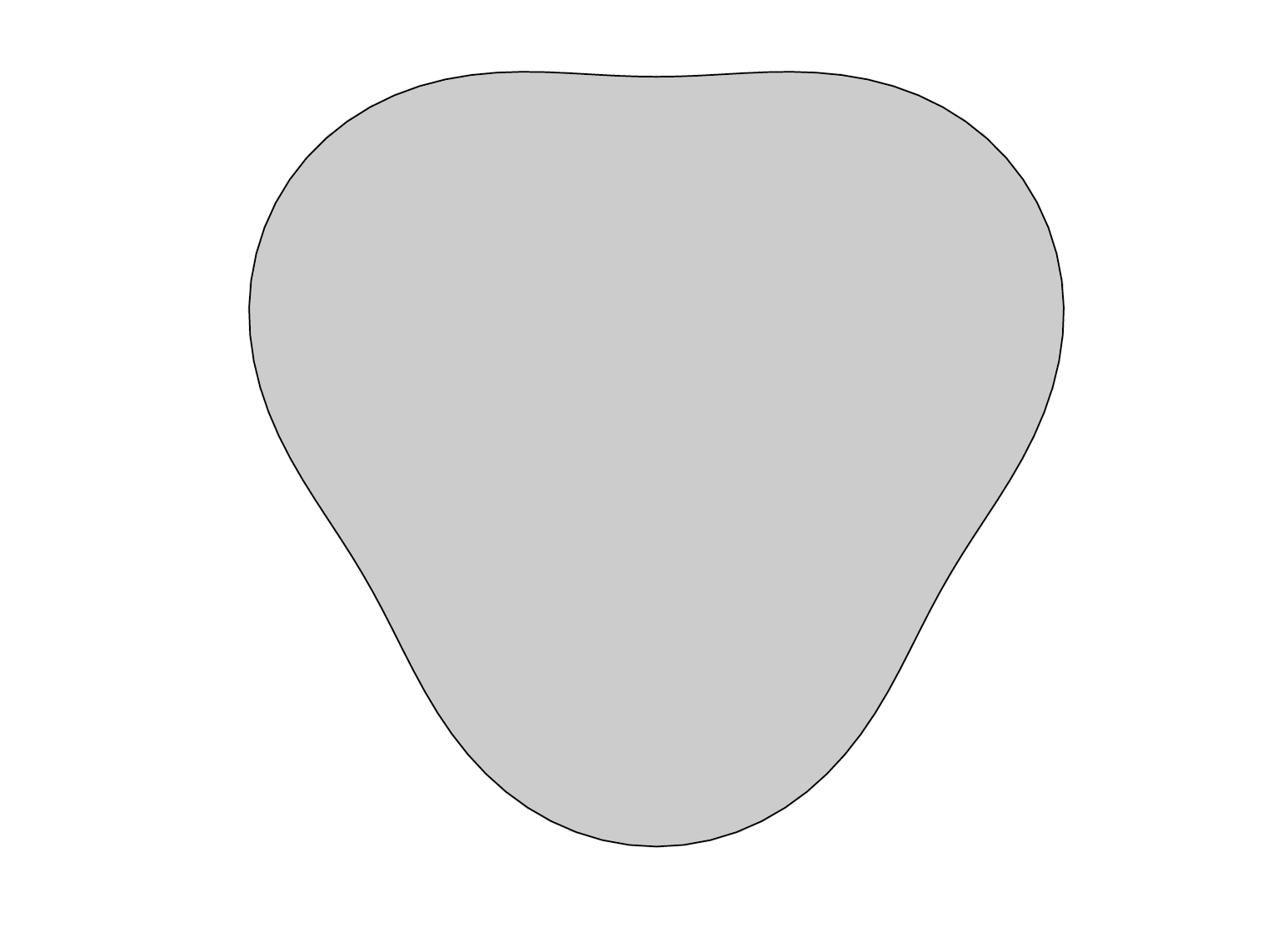}
\end{subfigure}\quad
\begin{subfigure}{}
\includegraphics[width=3.5cm]{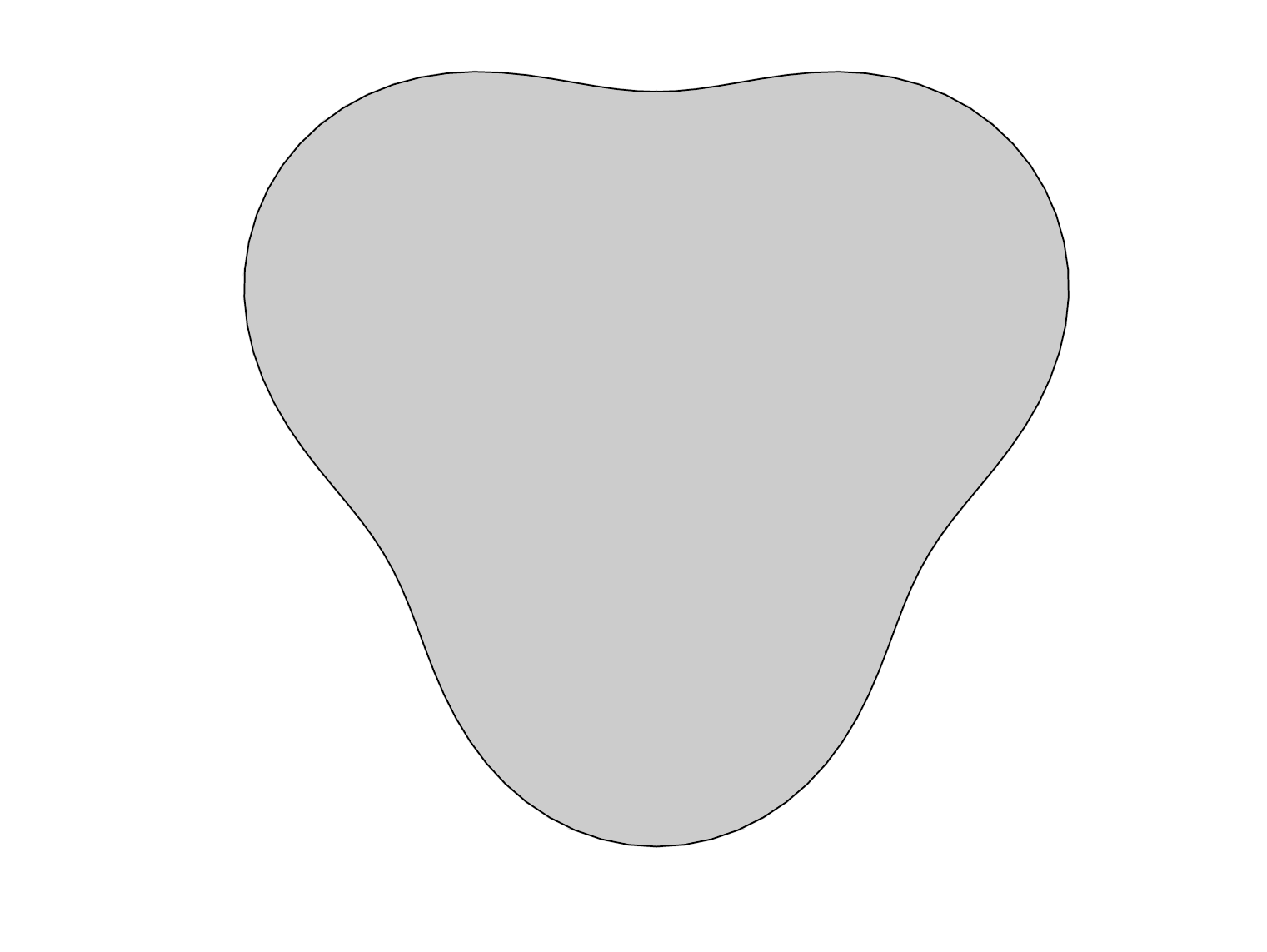}
\end{subfigure}\quad
\begin{subfigure}{}
\includegraphics[width=3.5cm]{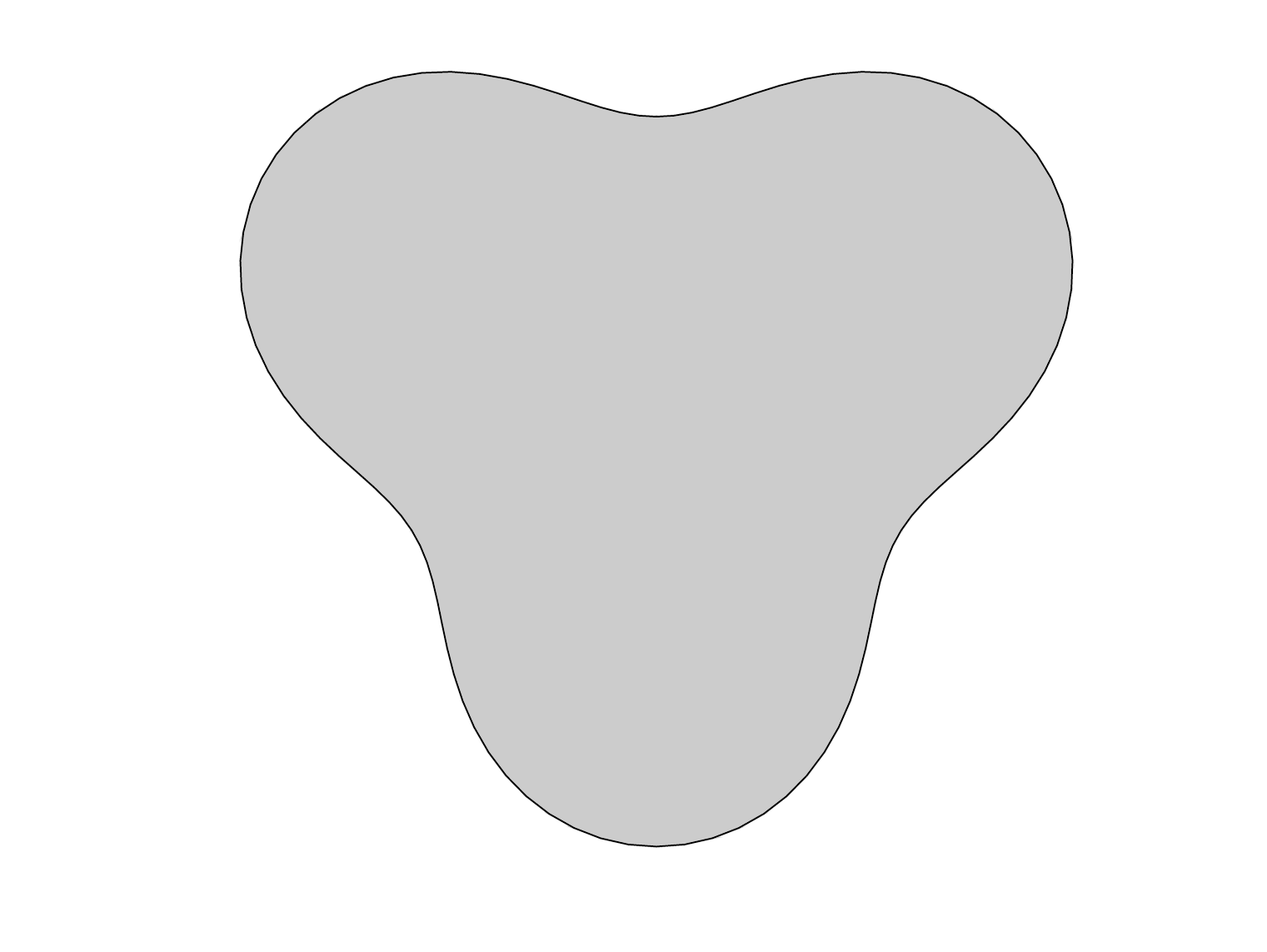}
\end{subfigure}\quad
\begin{subfigure}{}
\includegraphics[width=3.5cm]{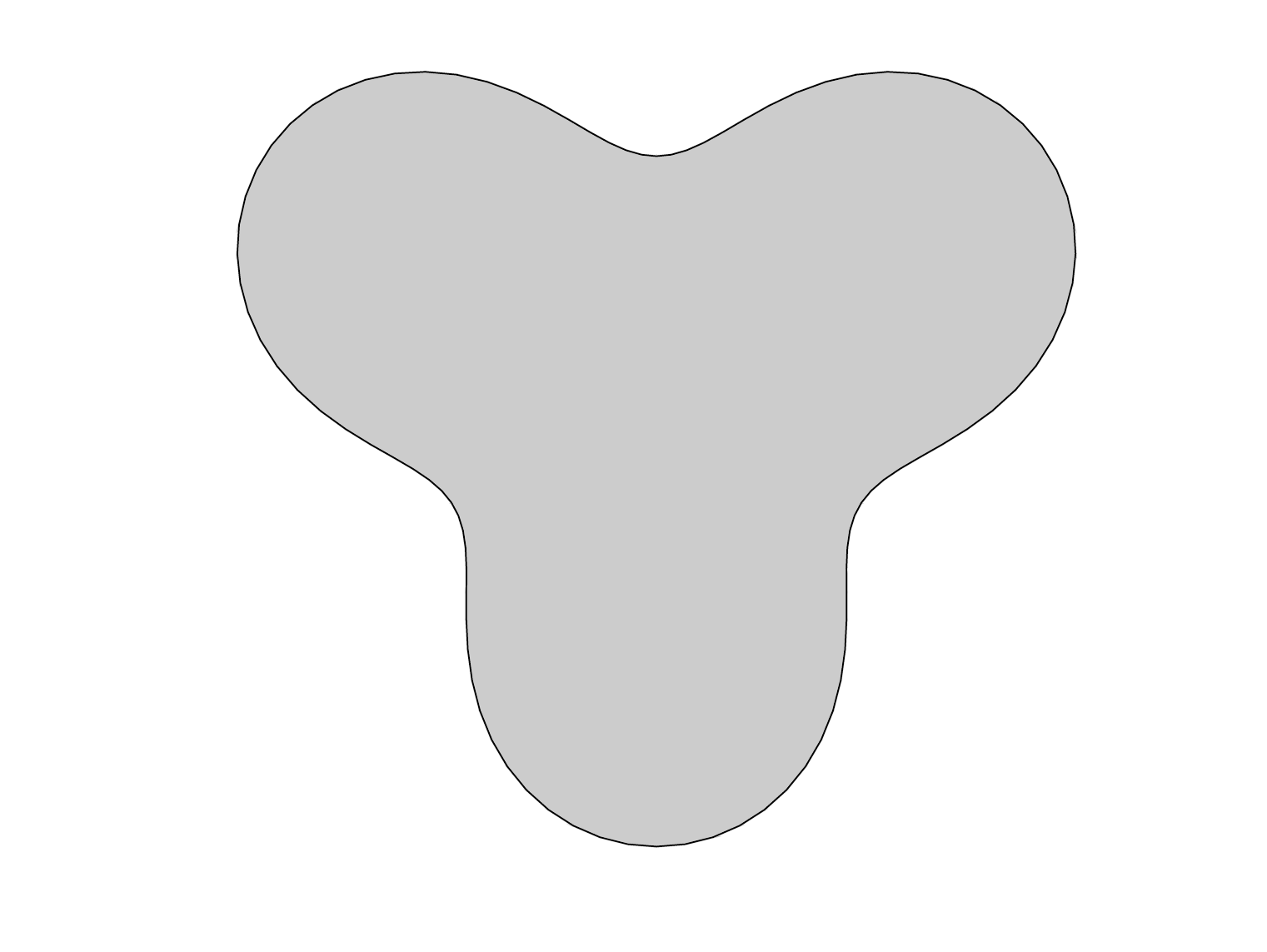}
\end{subfigure}\quad
\begin{subfigure}{}
\includegraphics[width=3.5cm]{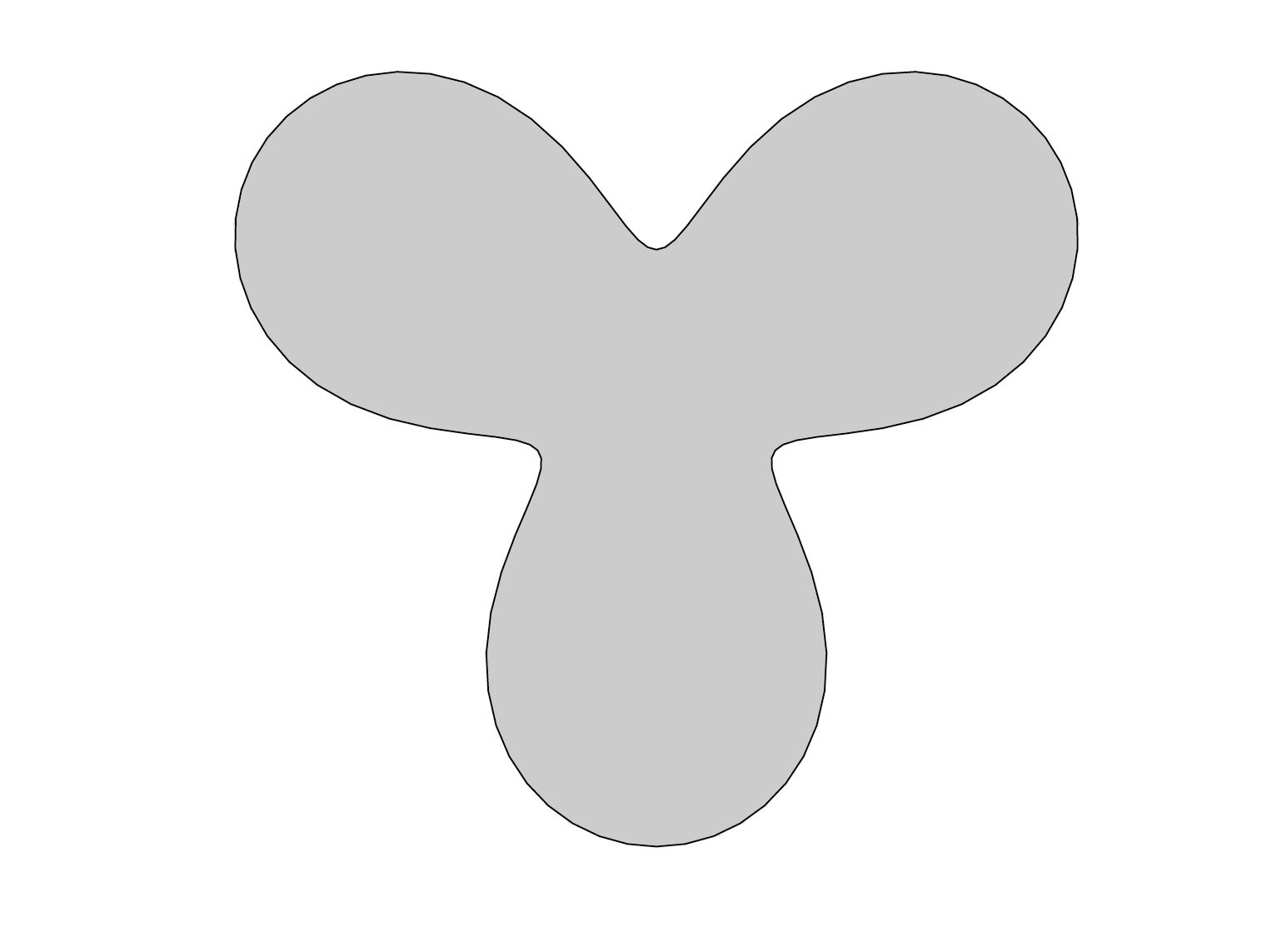}
\end{subfigure}\\
\hrule
\begin{subfigure}{}
\includegraphics[width=3.5cm]{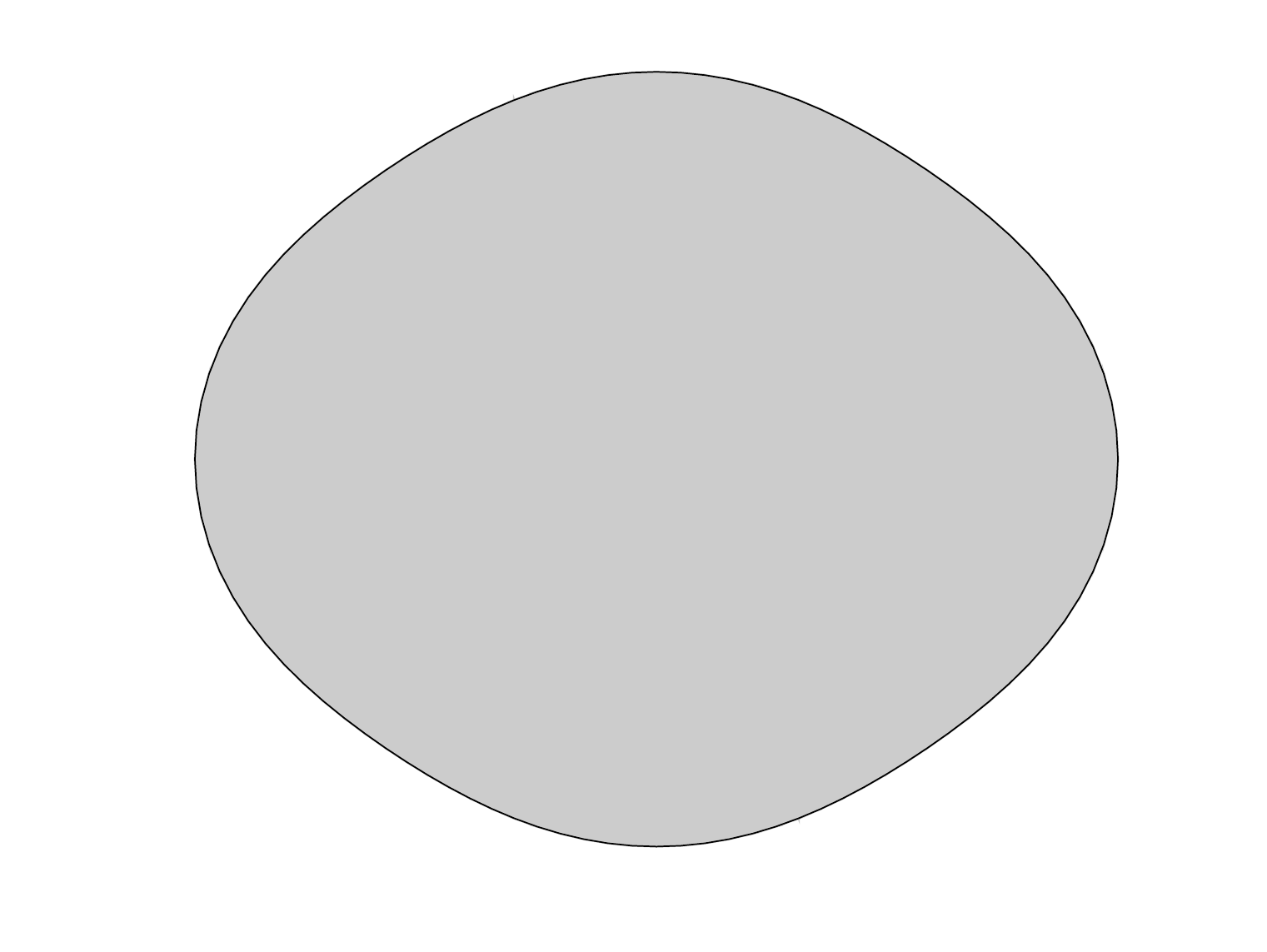}
\end{subfigure}\quad
\begin{subfigure}{}
\includegraphics[width=3.5cm]{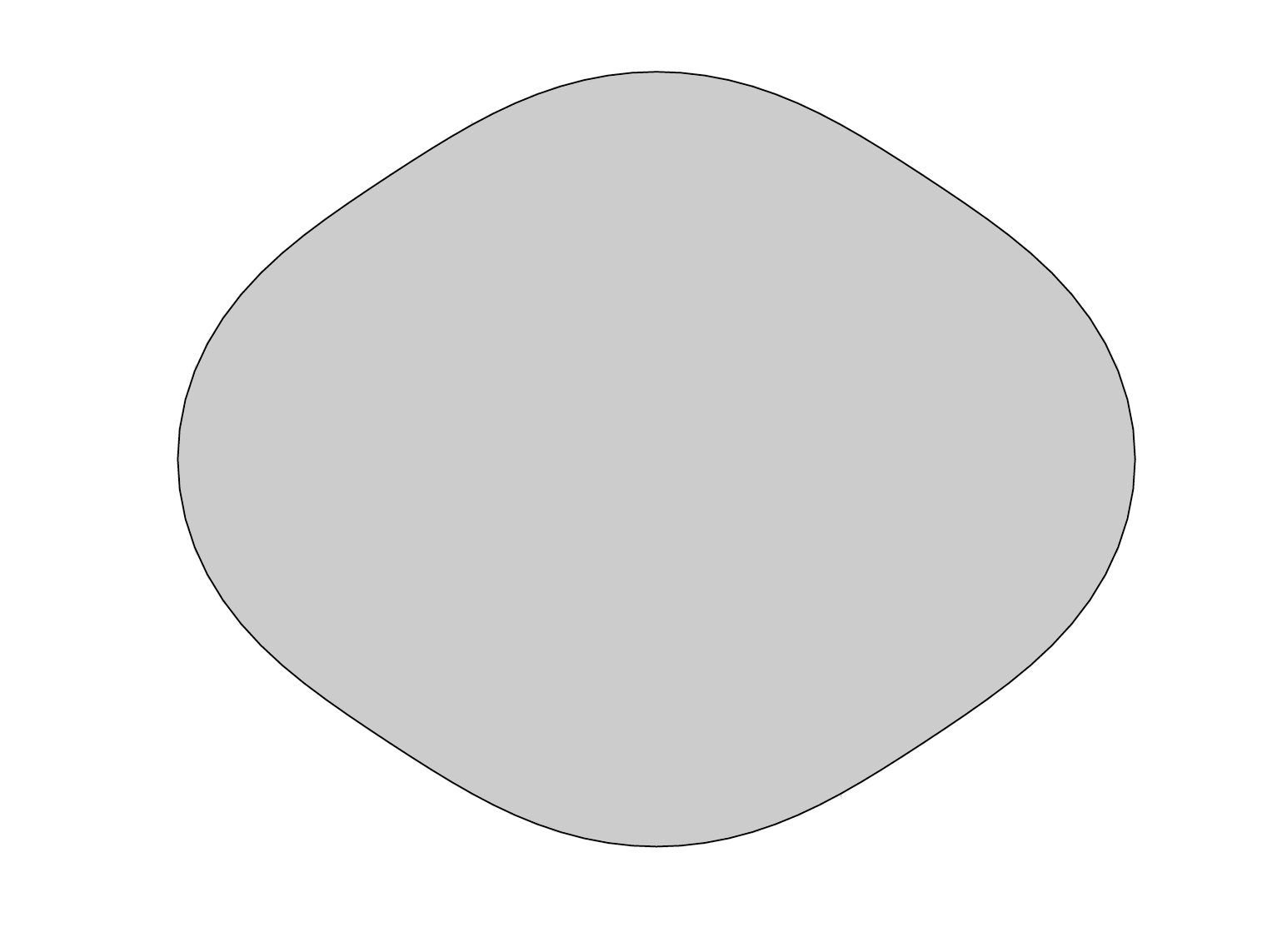}
\end{subfigure}\quad
\begin{subfigure}{}
\includegraphics[width=3.5cm]{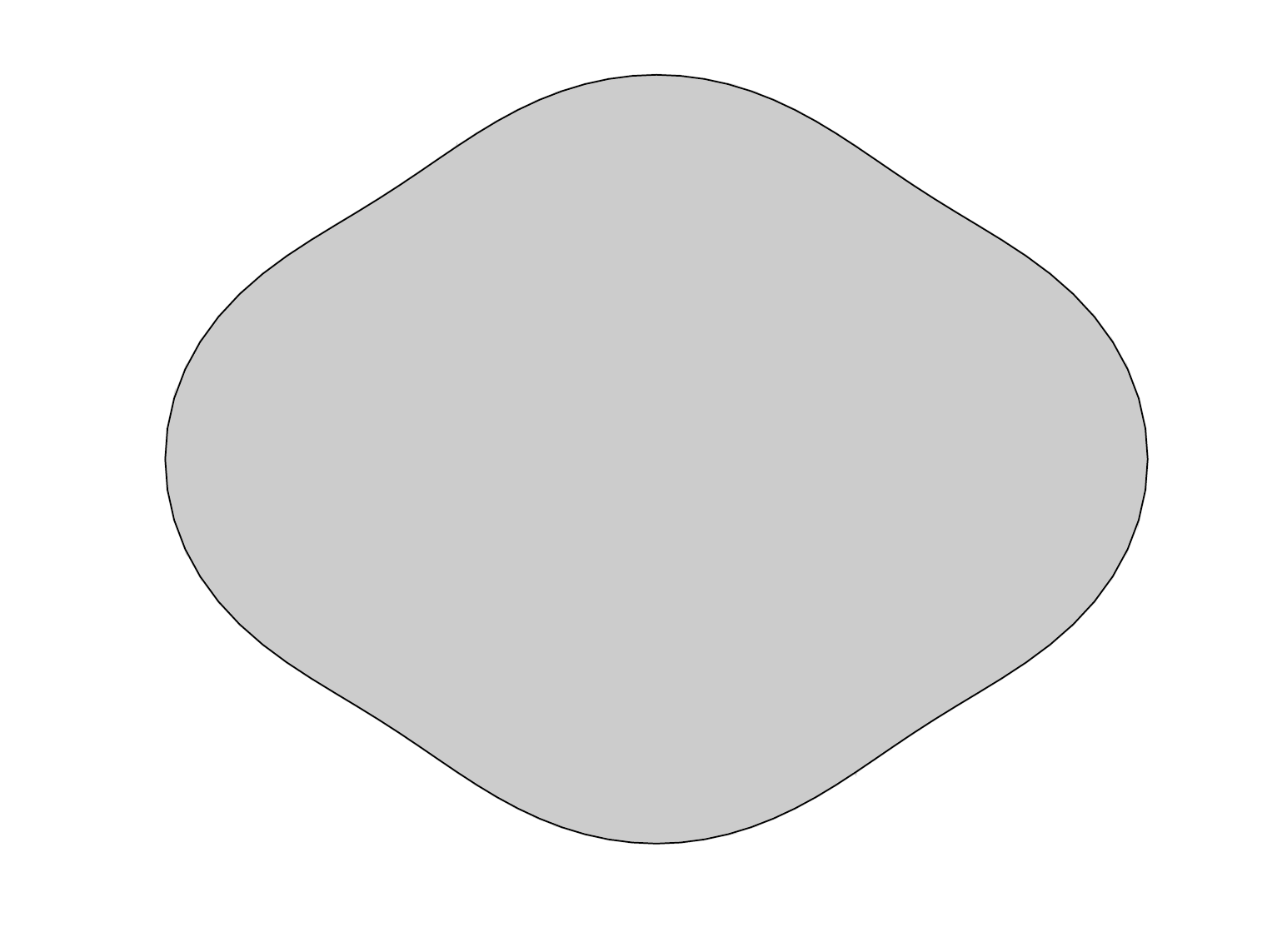}
\end{subfigure}\quad
\begin{subfigure}{}
\includegraphics[width=3.5cm]{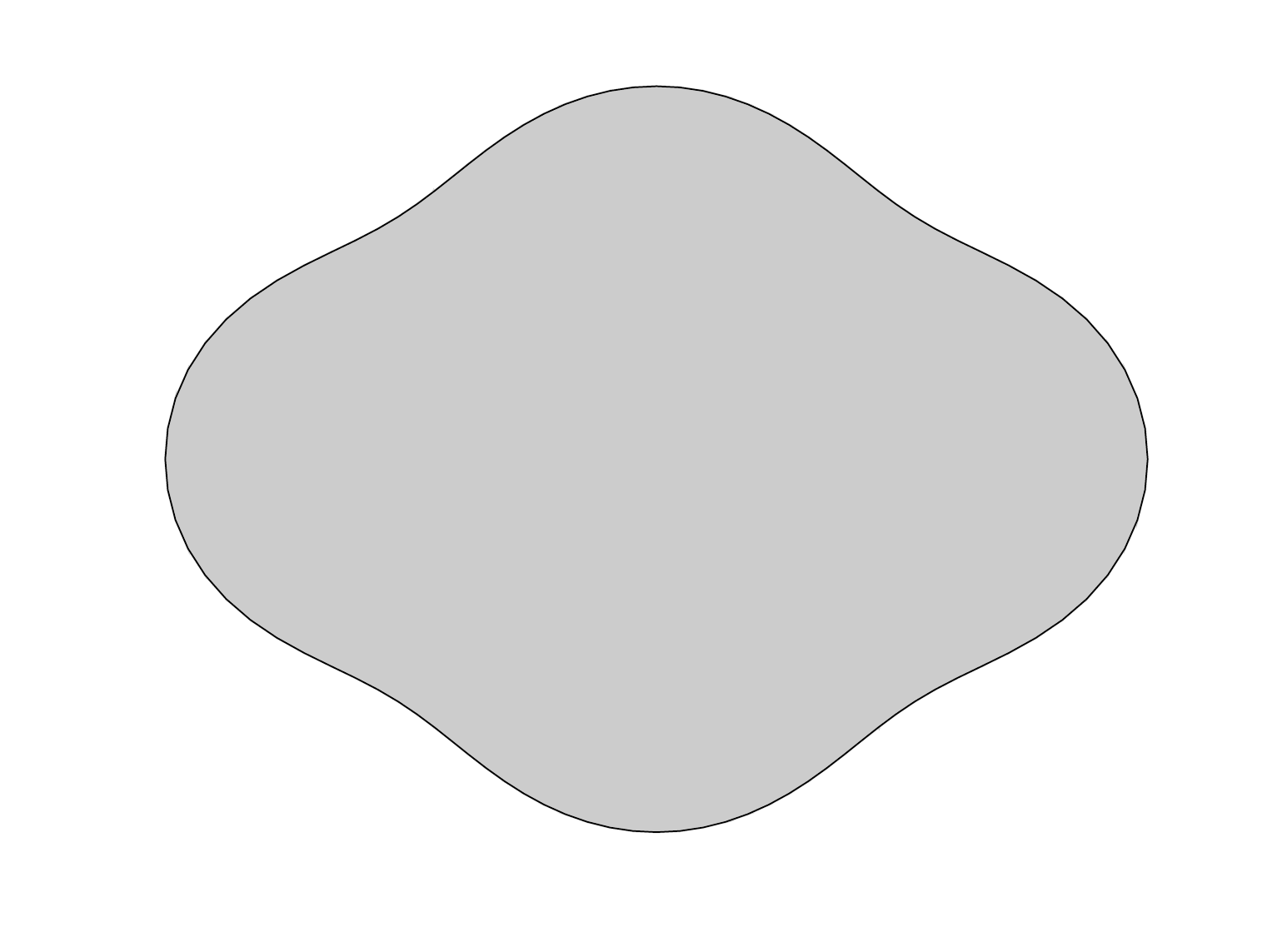}
\end{subfigure}\quad
\begin{subfigure}{}
\includegraphics[width=3.5cm]{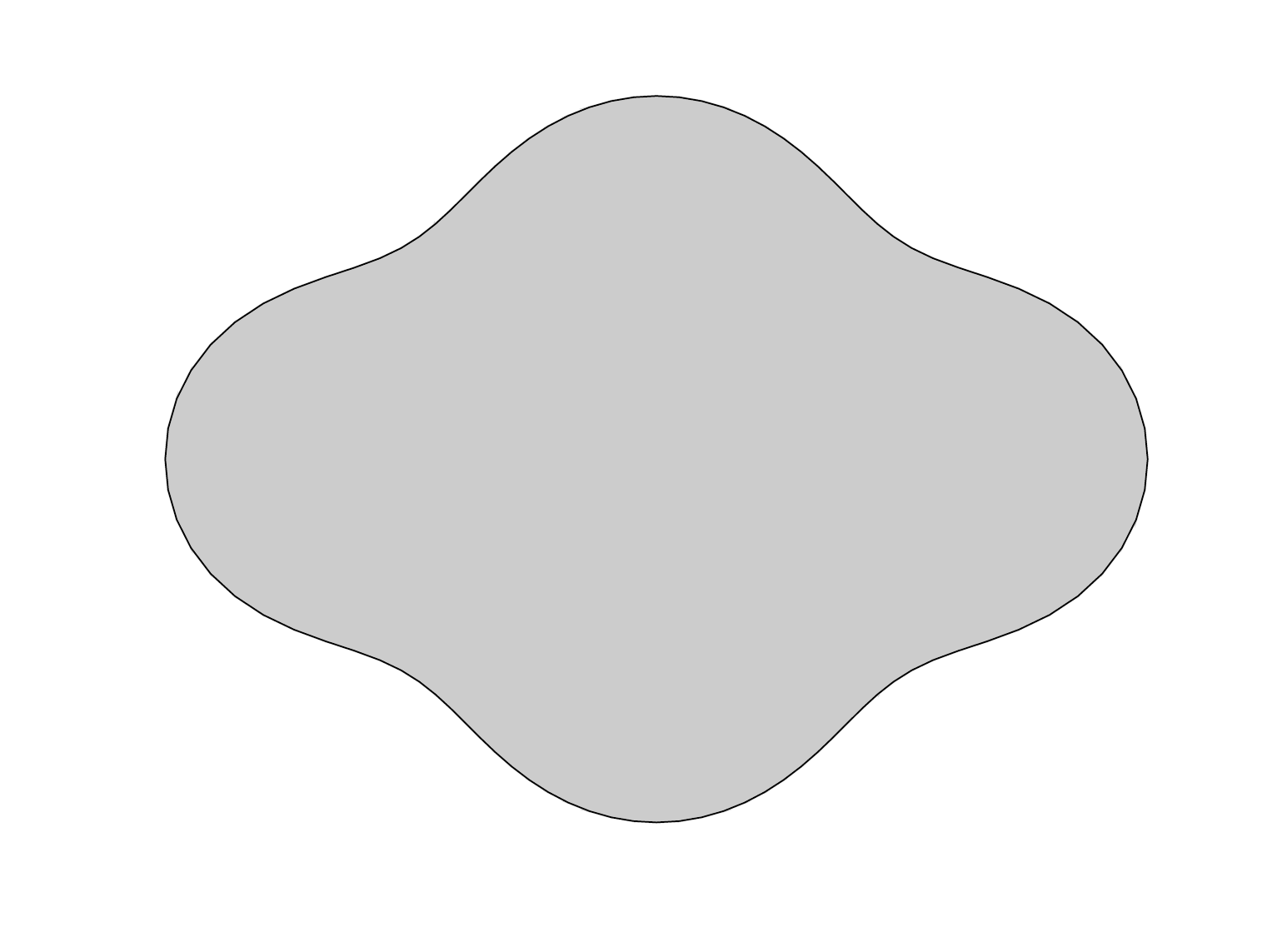}
\end{subfigure}\quad
\begin{subfigure}{}
\includegraphics[width=3.5cm]{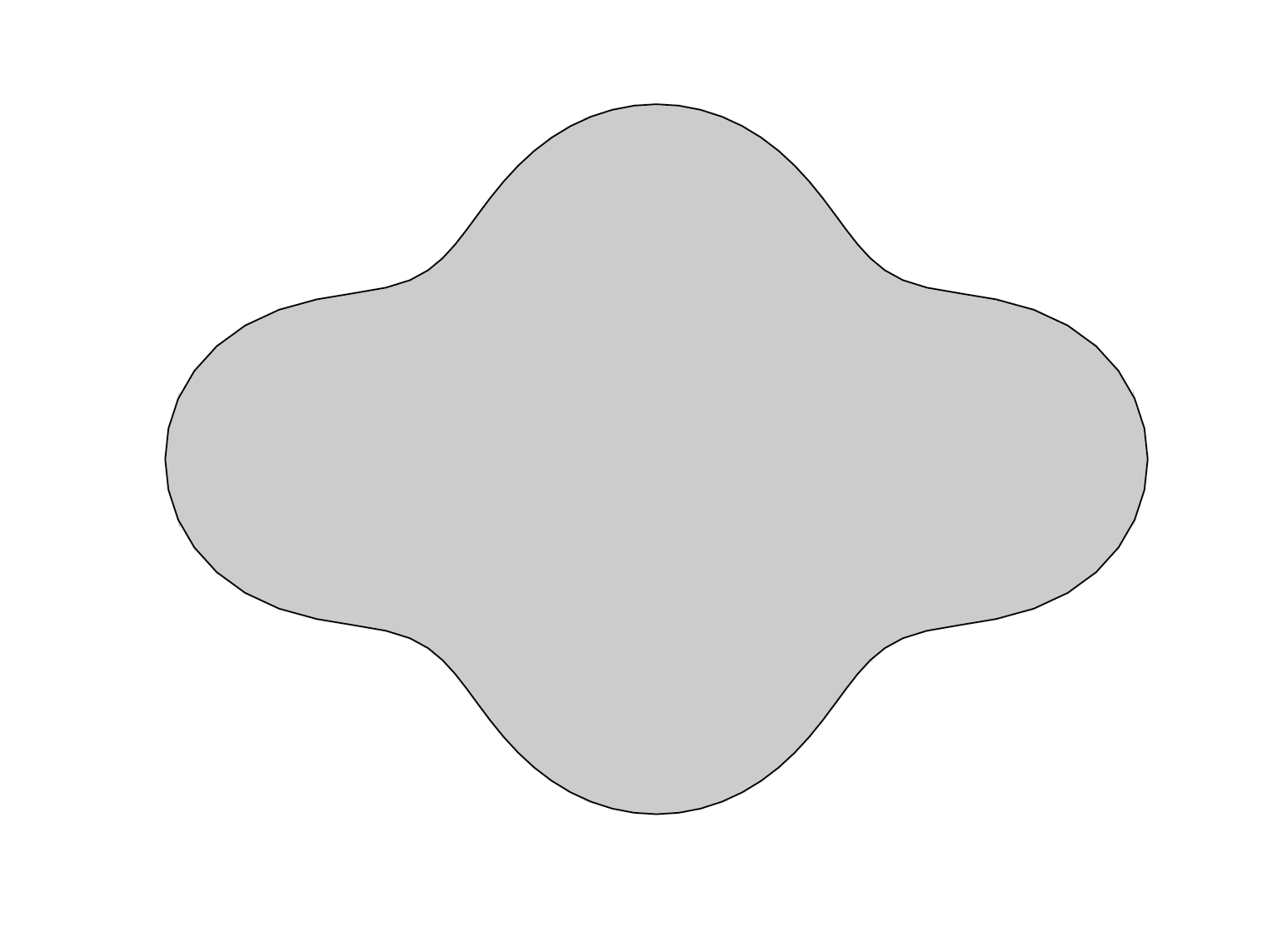}
\end{subfigure}\\
\caption{The influence of the parameter $c=1.75$, $2.00$, $2.25$, $2.50$, $2.75$, and $3.00$ for $m=3$ (first and second row) and for $m=4$ (third and fourth row).}
\label{fig:kle3}
\end{figure}
Additionally, one can see that we are almost able to obtain a possible shape of the maximizer for the third and fourth INE. 
To add more flexibility, we introduce the additional parameter $\alpha$. The modified equipotentials are given in the form
\begin{eqnarray}
\sum_{i=1}^m \frac{1}{\|x-P_i\|^{2\alpha}}=c
\label{impl:kle}
\end{eqnarray}
We introduce the two in front of the parameter $\alpha$ in order to avoid the computation of the square root in the norm definition.
In Figure \ref{fig:kle4} we show the influence of the parameter $\alpha$ fixing $c=2$. As one can see, we have enough flexibility to obtain very good approximations for a possible shape maximizer for the third and fourth INE.
\begin{figure}
\centering
\begin{subfigure}{}
\includegraphics[width=3.5cm]{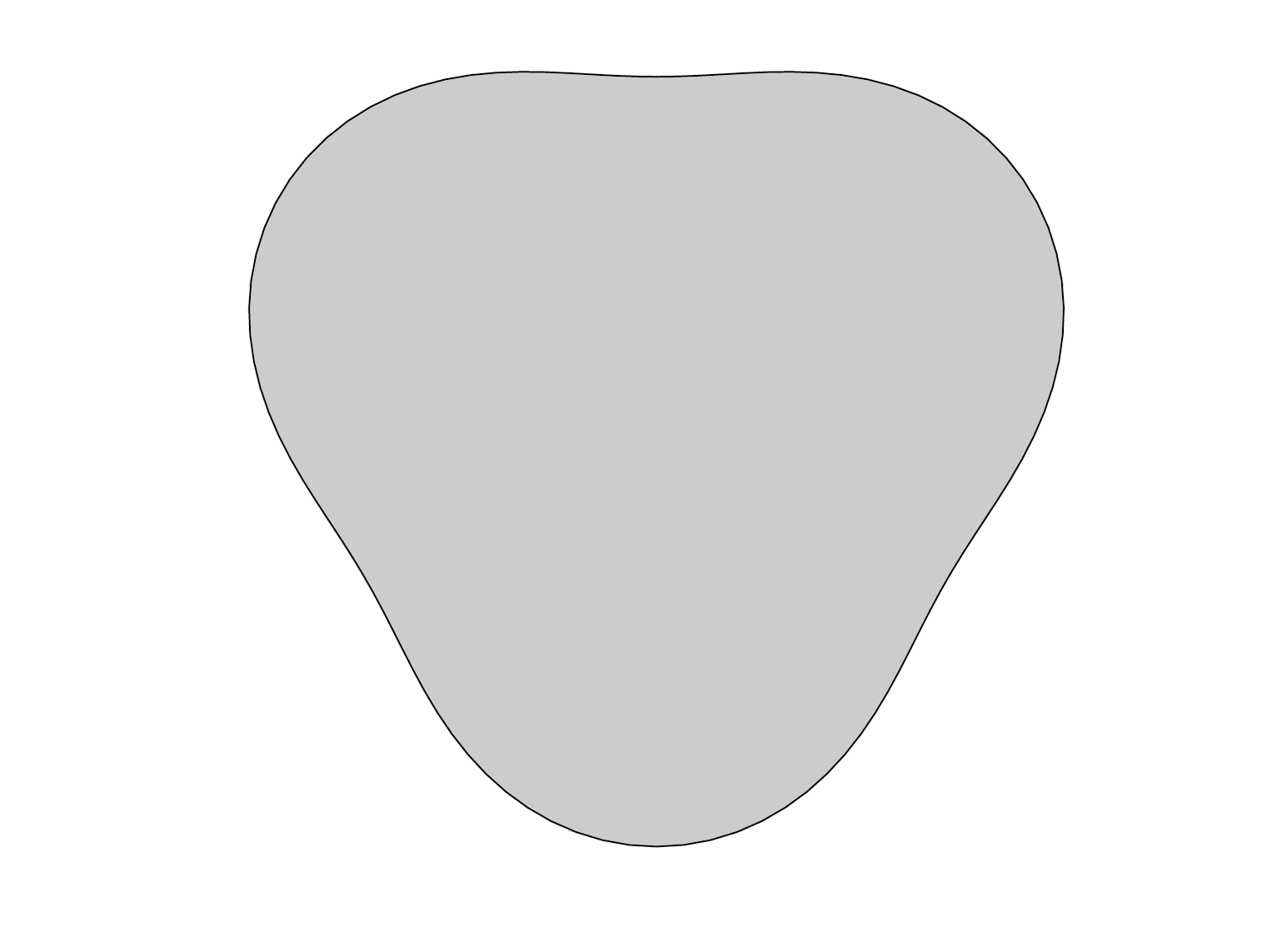}
\end{subfigure}\quad
\begin{subfigure}{}
\includegraphics[width=3.5cm]{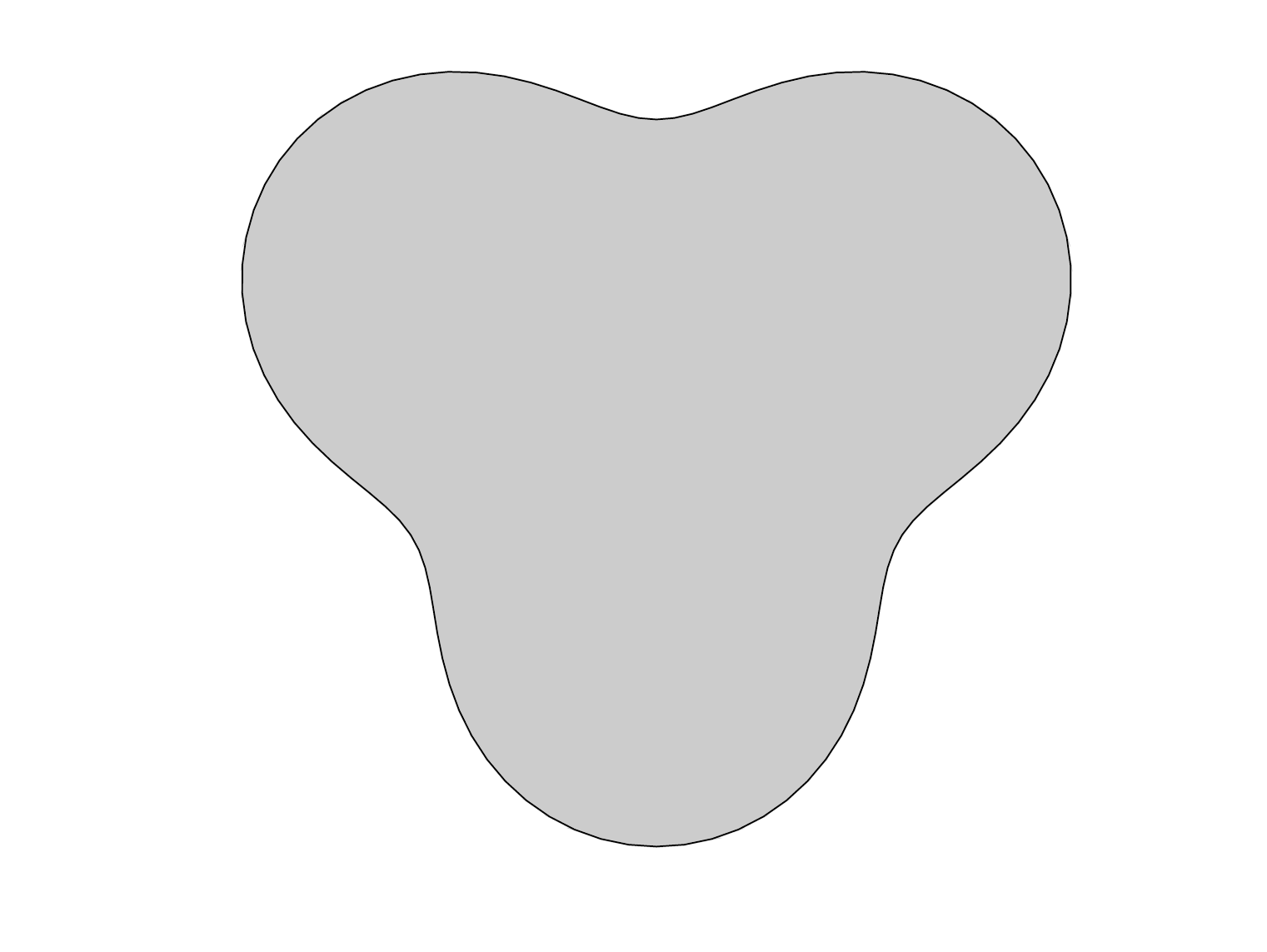}
\end{subfigure}\quad
\begin{subfigure}{}
\includegraphics[width=3.5cm]{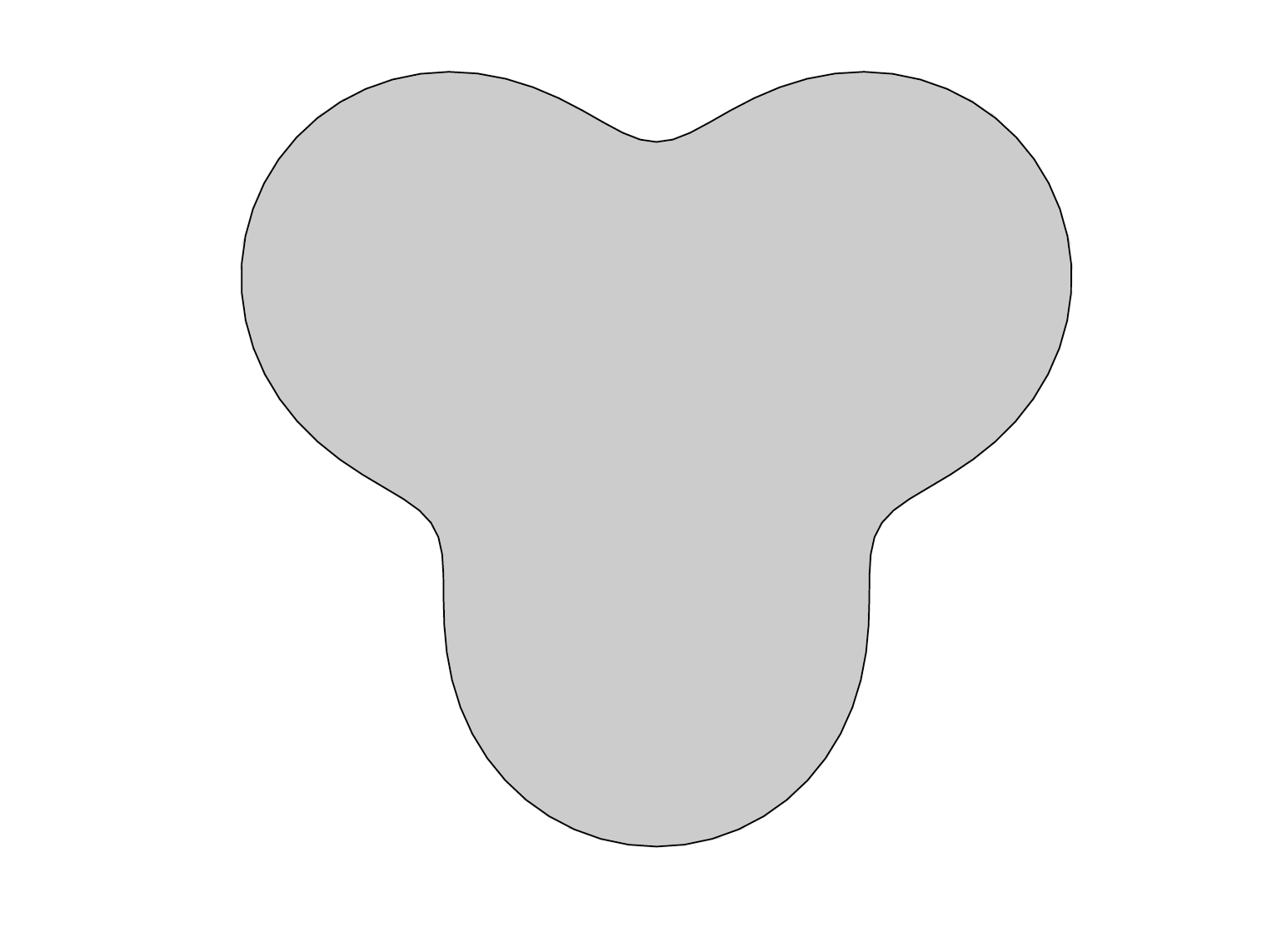}
\end{subfigure}\quad
\begin{subfigure}{}
\includegraphics[width=3.5cm]{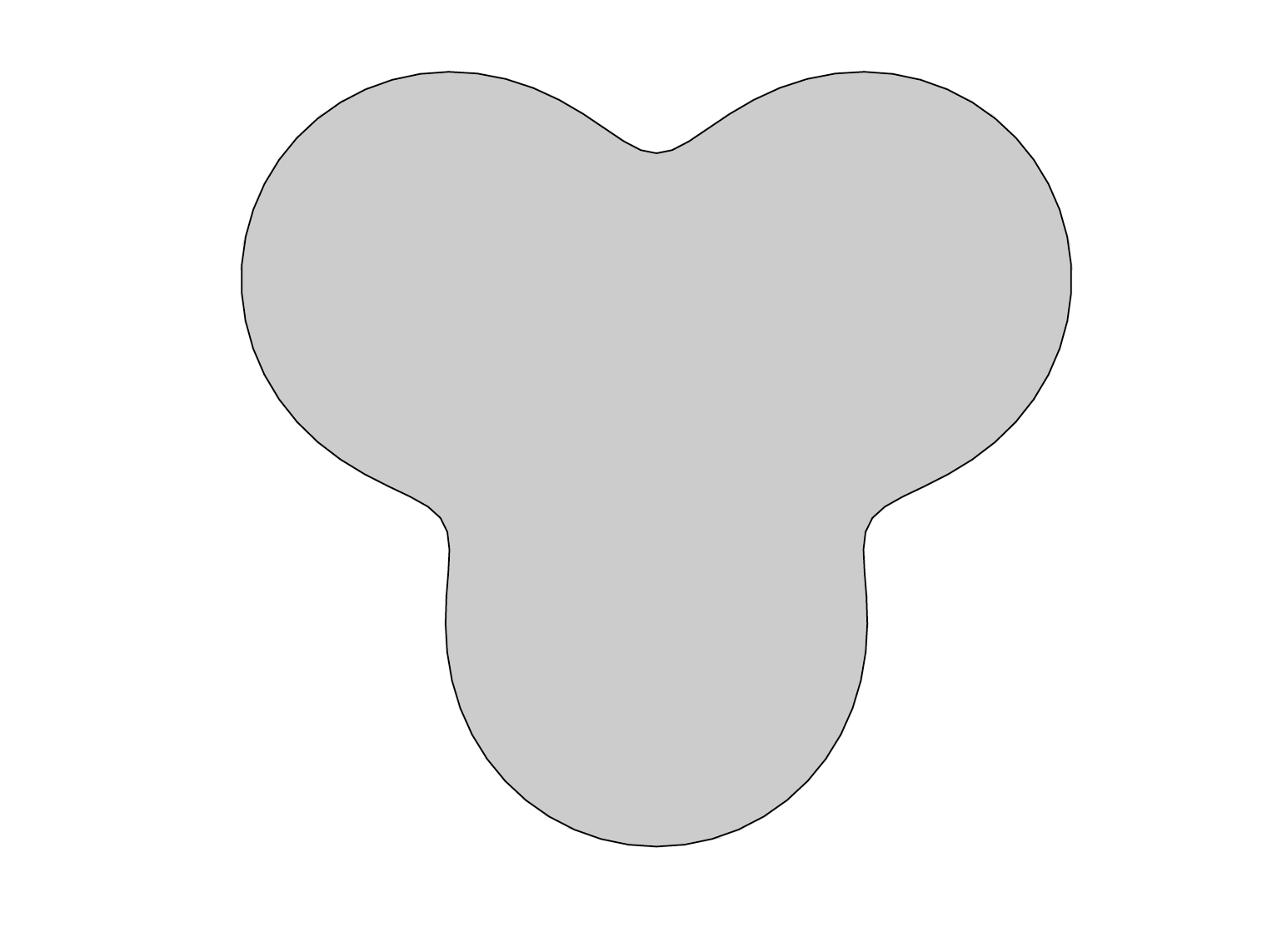}
\end{subfigure}\quad
\begin{subfigure}{}
\includegraphics[width=3.5cm]{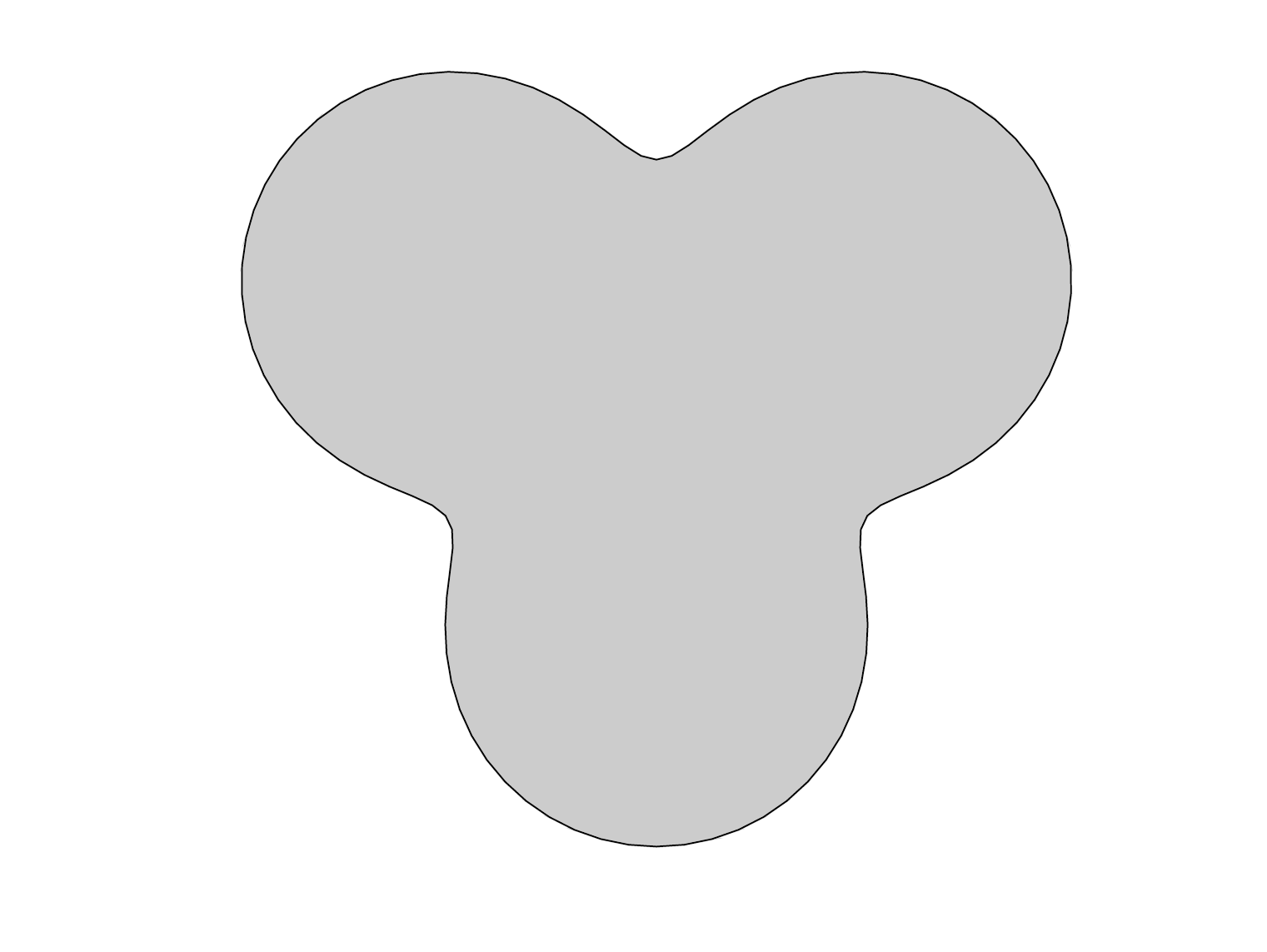}
\end{subfigure}\quad
\begin{subfigure}{}
\includegraphics[width=3.5cm]{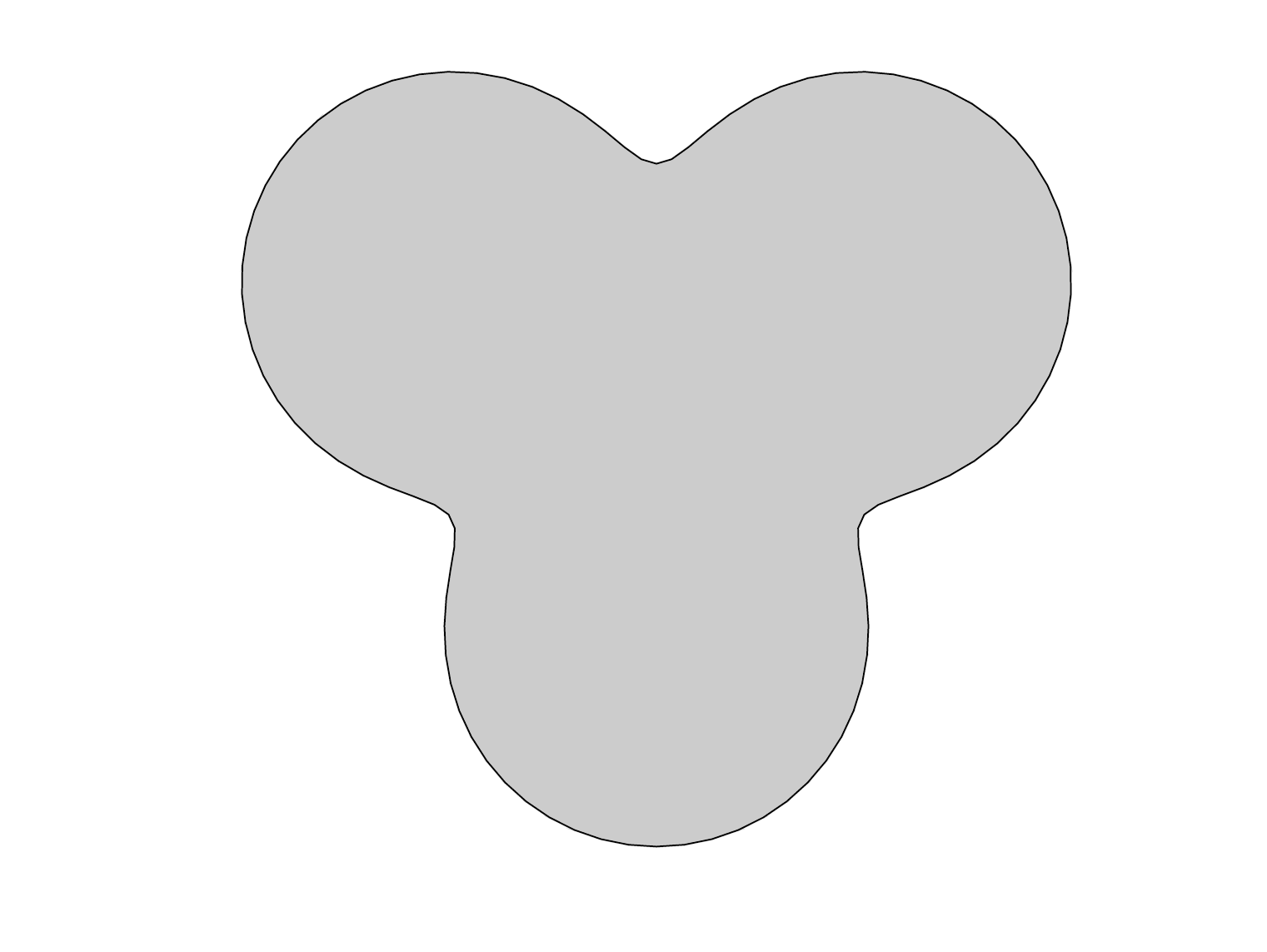}
\end{subfigure}\\
\hrule
\begin{subfigure}{}
\includegraphics[width=3.5cm]{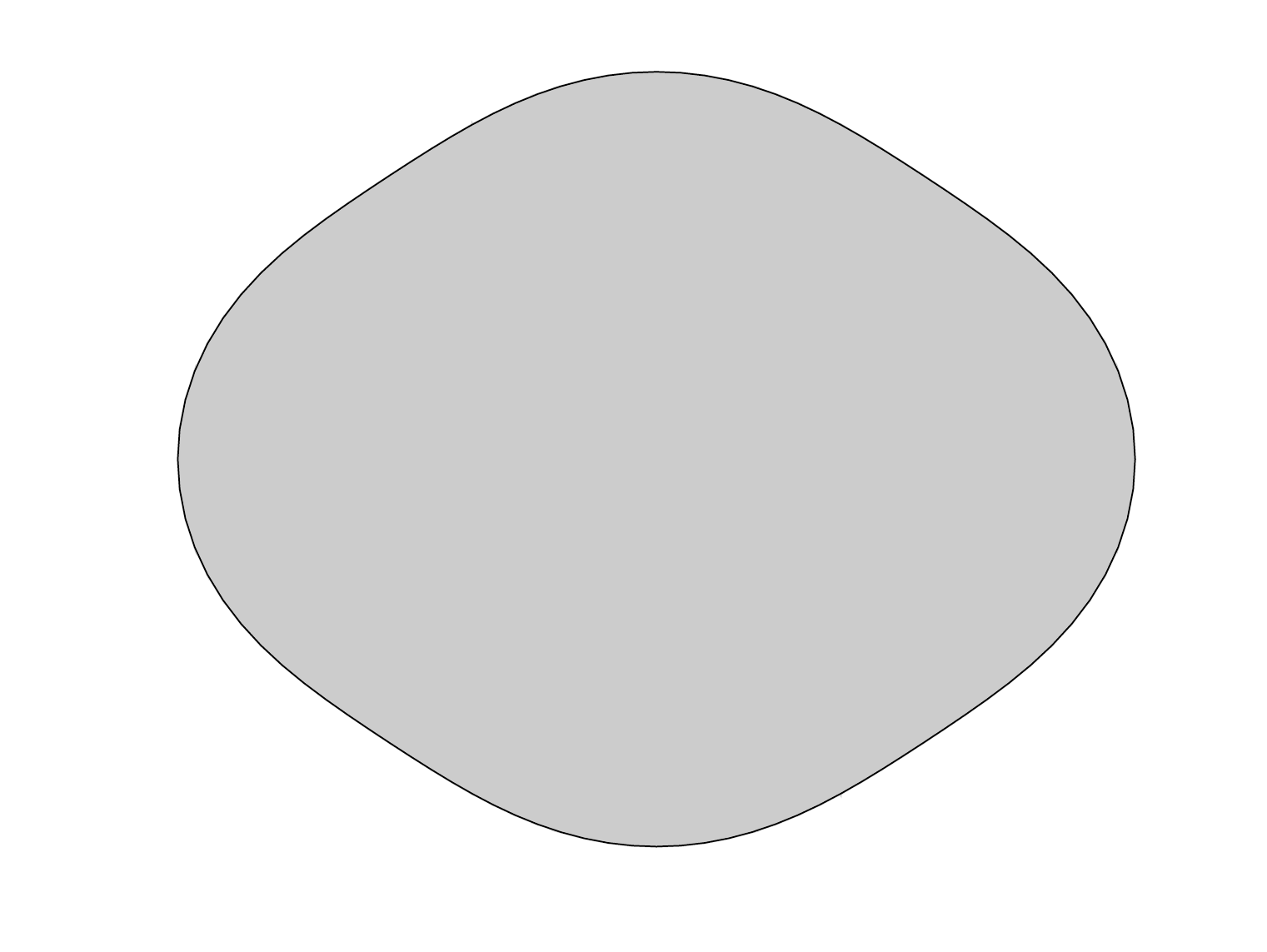}
\end{subfigure}\quad
\begin{subfigure}{}
\includegraphics[width=3.5cm]{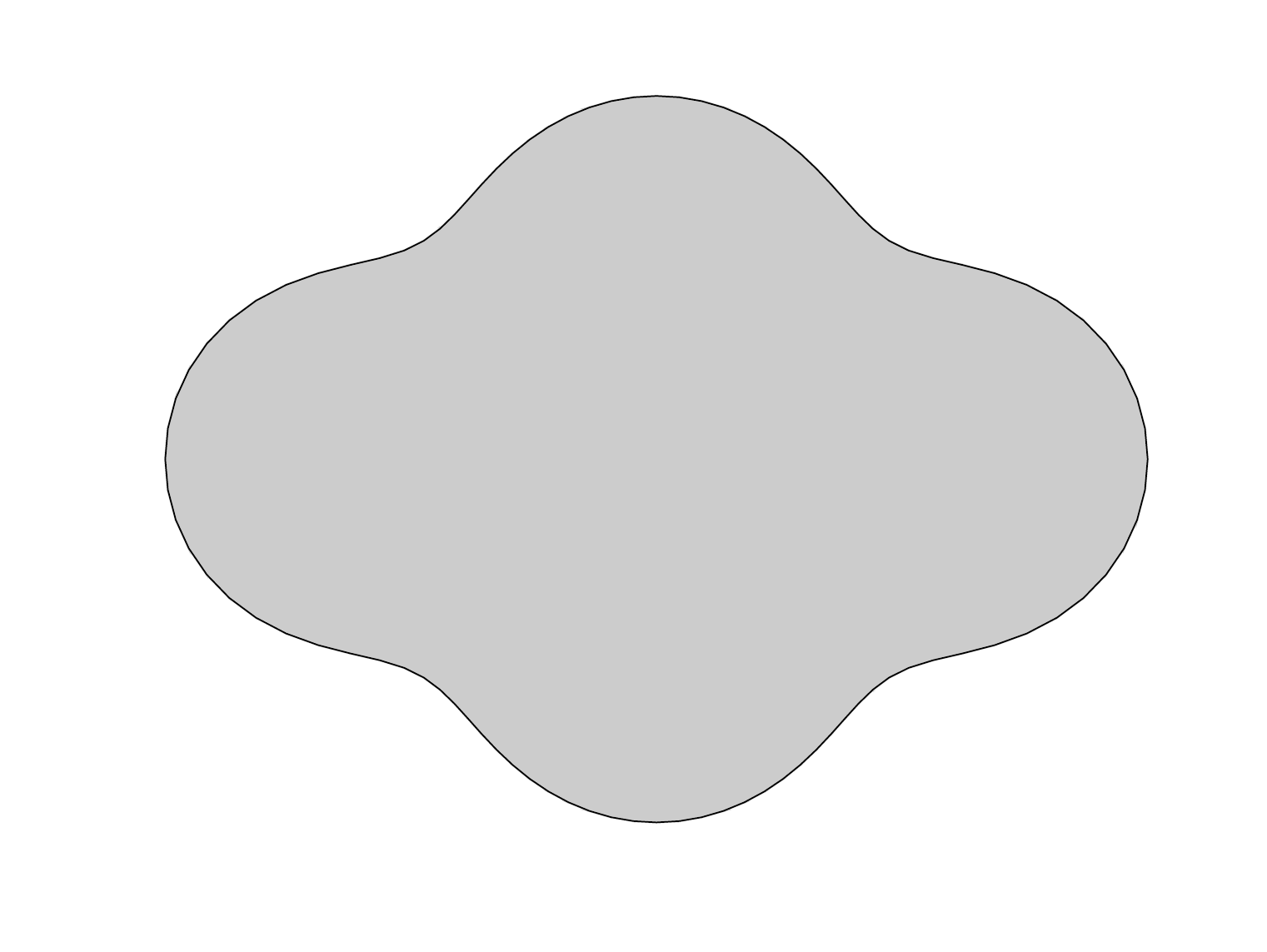}
\end{subfigure}\quad
\begin{subfigure}{}
\includegraphics[width=3.5cm]{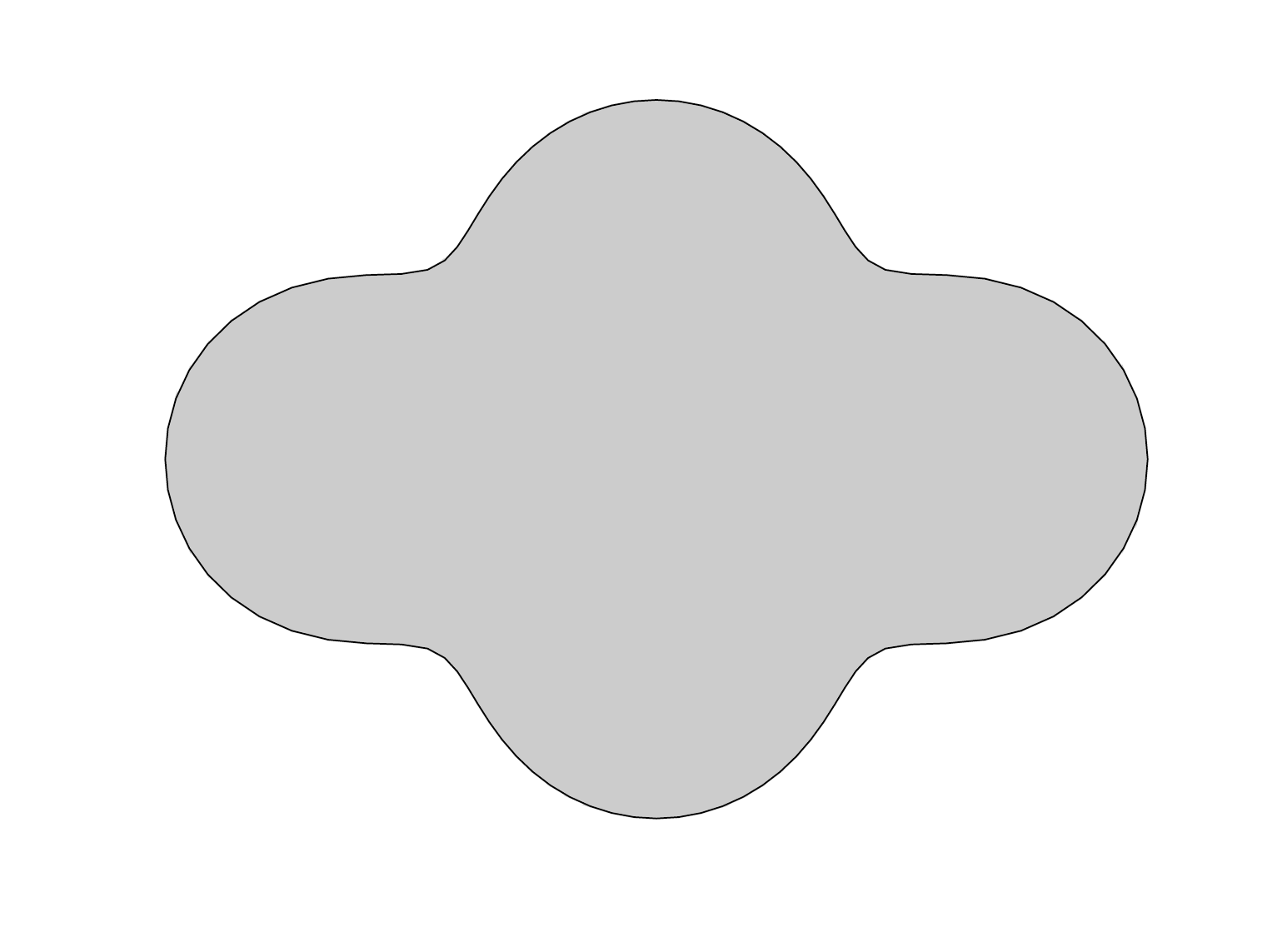}
\end{subfigure}\quad
\begin{subfigure}{}
\includegraphics[width=3.5cm]{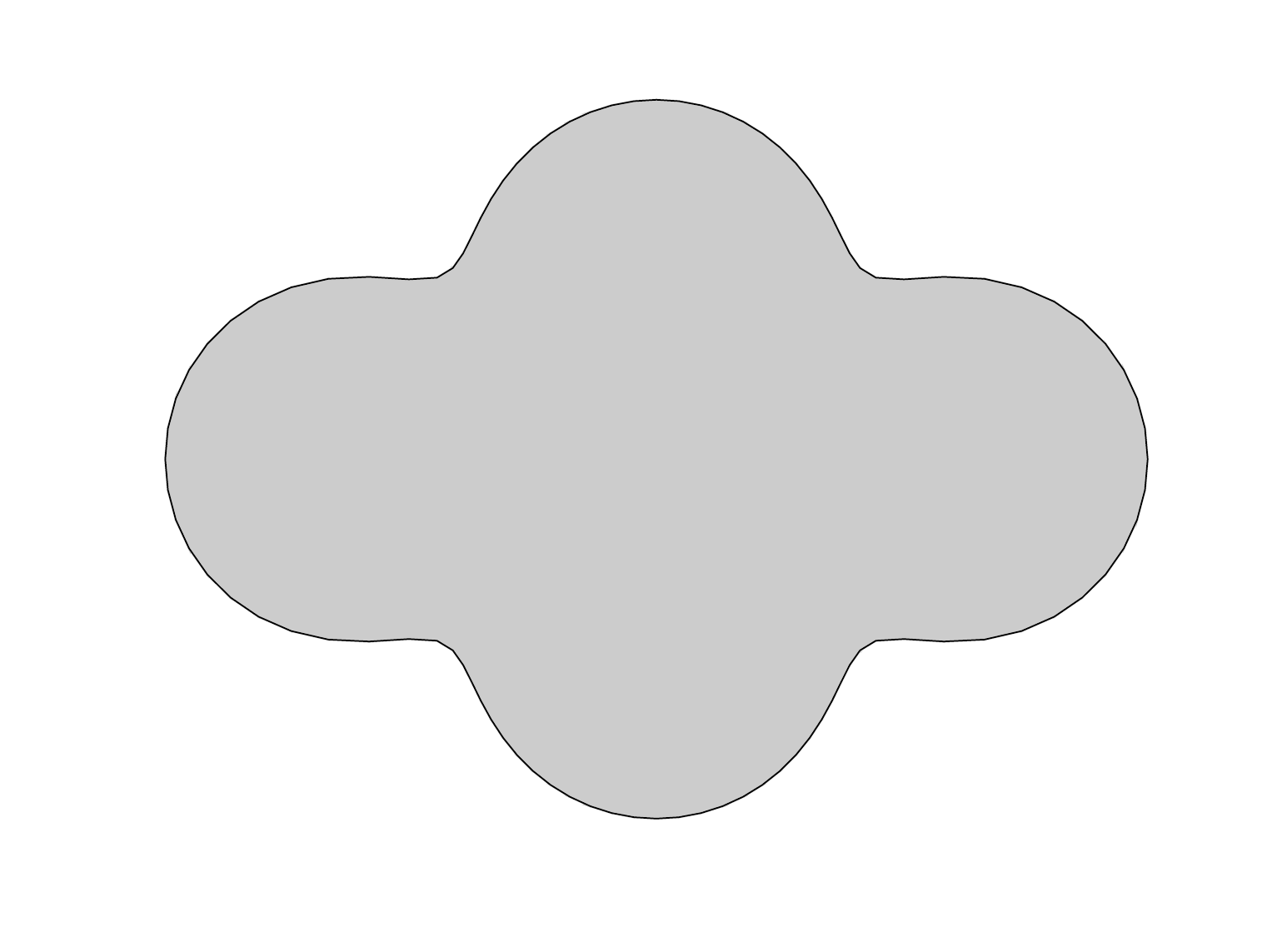}
\end{subfigure}\quad
\begin{subfigure}{}
\includegraphics[width=3.5cm]{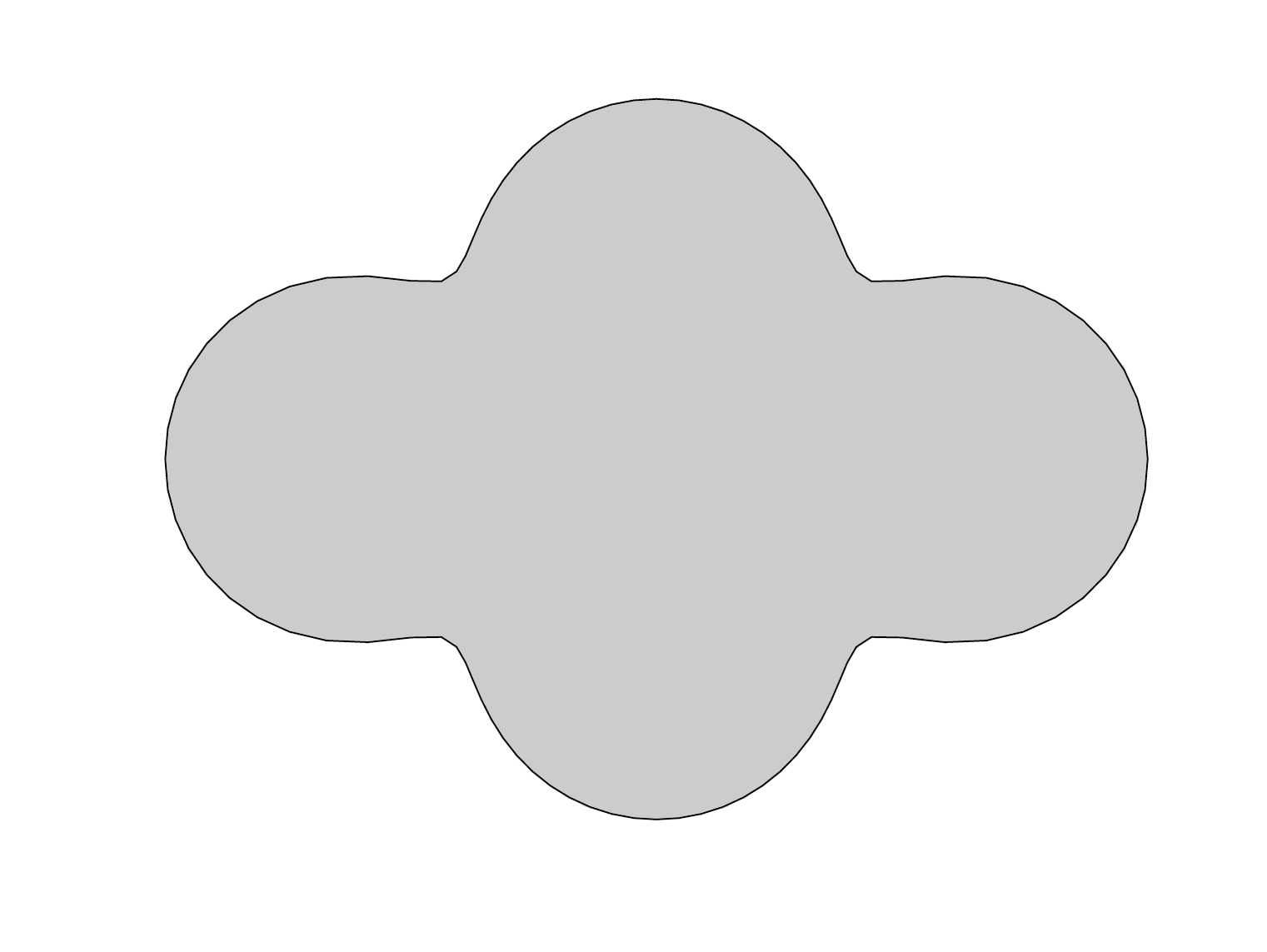}
\end{subfigure}\quad
\begin{subfigure}{}
\includegraphics[width=3.5cm]{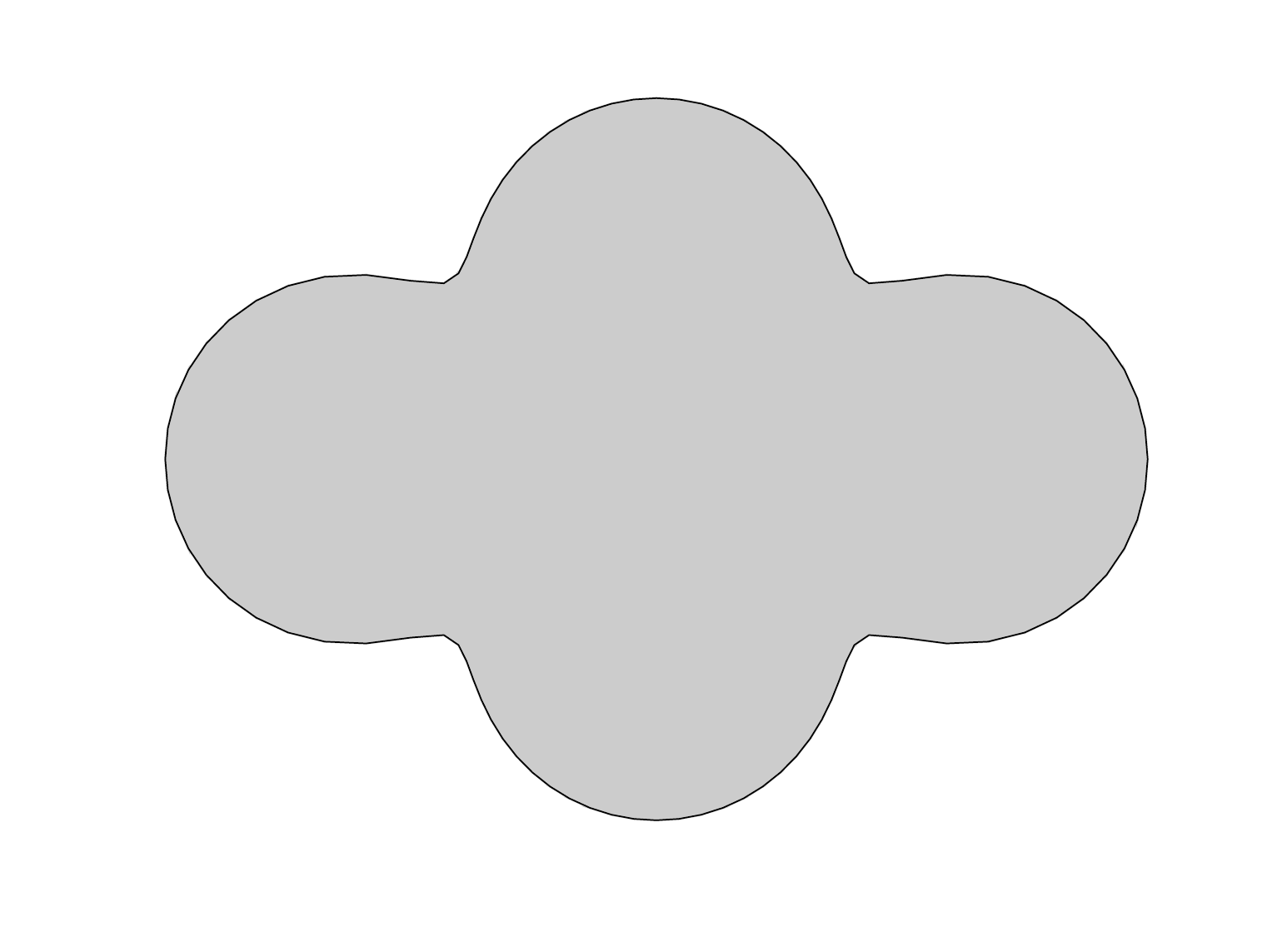}
\end{subfigure}\\
\caption{The influence of the parameter $\alpha=0.5$, $1.0$, $1.5$, $2.0$, $2.5$, and $3.0$ with fixed $c=2$ for $m=3$ (first and second row) and for $m=4$ (third and fourth row).}
\label{fig:kle4}
\end{figure}
Thus, we have seen the influence of the parameters $\alpha$ and $c$. We shortly explain how to generate $n$ points on the boundary for the given parameters $\alpha$ and $c$.
This is done as follows. First, the equation (\ref{impl:kle}) is rewritten in polar coordinates. Then, $n+1$ equidistant angle $\phi$ in the interval $[0,2\pi]$ are generated. 
Next, for each angle $\phi_i$ the implicit equation is solved for the unique $r_i$ via a root finding algorithm. 
Finally, the points given in polar coordinates $(r_i,\phi_i)$, $i=1,\ldots,n+1$ are transformed back to rectangular coordinates $(x_i,y_i)=(r_i\cos(\phi_i),r_i\sin(\phi_i))$, $i=1,\ldots,n+1$. 
Hence, we obtain $n$ different points on the boundary of the scatterer (the $(n+1)$-th point is the same as the first point by construction). 
Those $n$ points can now be used in the boundary element collocation method. 

In order to calculate the value $\lambda_k\cdotp A$, we need to numerically approximate the area enclosed by the given implicit curve (see (\ref{impl:kle})).
That is, we have $n$ points distributed on the boundary $\partial \Omega$. With these points and the approximation via quadratic interpolation, the domain $\widetilde{\Omega}$ with the boundary 
$\widetilde{\partial \Omega}$ is defined. To approximate the area of this region, we compute the area of the non-self intersecting polygon spanned by choosing $p\gg n$ points including an additional point 
(the first point is the additional $(p+1)$-th point). 
The approximate area is given by
\[A\approx A_{\widetilde{\Omega}}=\frac{1}{2}\left|\sum_{i=1}^p(x_i-x_{i+1})(y_i+y_{i+1})\right|\]
which is an easy consequence of the formula (\cite[4.6.1, p. 206]{Zw12})
\[\frac{1}{2}\left|\left|\begin{matrix}x_1 & x_2\\ y_1 & y_2\end{matrix}\right|+
                   \left|\begin{matrix}x_2 & x_3\\ y_2 & y_3\end{matrix}\right|+\ldots+
                   \left|\begin{matrix}x_p & x_1\\ y_p & y_1\end{matrix}\right|\right|\,.\]

The exterior normals on the boundary given implicitly by (\ref{impl:kle}) are given by $\nu=\tilde{\nu}/\|\tilde{\nu}\|$ with 
 \begin{eqnarray*}
 \tilde{\nu}=-2\alpha \sum_{i=1}^m \frac{(x-P_i)}{\|x-P_i\|^{2(\alpha+1)}}\,.
 \end{eqnarray*}

Now, we have everything together in order to optimize with respect to the two parameter $c$ and $\alpha$. First, we consider the third INE. 
The reference value given by Antunes \& Oudet is given by $32.90$ using 37 unknown coefficients. The third eigenvalue has multiplicity three. If we fix $\alpha=3/2$, then the optimization with respect to $c$ yields the result $c=1.8416$ with $32.8929$, $32.8929$, $32.8929$ for the third, fourth, and fifth, respectively. 
As we observe, the reported numbers are more accurate. If we fix $\alpha=2$, then we obtain $c=1.6921$ with $32.9018$, $32.9018$, $32.9018$ which improves the result slightly compared to the value $32.90$. 
But remember that we have only one unknown describing the boundary. If we choose $\alpha=5/2$, then we have $c=1.6112$ with $32.8970$, $32.8970$, $32.8970$. If we optimize with respect to both parameters yields $\alpha=2.0171$ and $c=1.6883$ with $32.9018$, $32.9018$, $32.9018$. The situation slightly changes for the optimization of the fourth eigenvalue. The reference value of Antunes \& Oudet is given by $43.86$ with multiplicity three using 33 unknown coefficients. 
If we use $\alpha=2$, we obtain  $c=2.0571$ with $43.6968$, $43.6968$, $44.2247$. Using $\alpha=5/2$ gives $c=2.0794$ with $43.8586$, $43.8586$, $43.8935$ which is close to the value of Antunes \& Oudet, 
but we have room for more considering the last eigenvalue. Fixing $\alpha=3$ yields $c=2.0875$ with $43.7822$, $43.7822$, $44.0634$. Optimizing with respect to the two parameters $\alpha$ and $c$ gives $\alpha=2.5426$ and $c=2.0845$ with $43.8694$, $43.8694$, $43.8694$. This is a much better result. 
In Figure \ref{fig:kle5} we show the three eigenfunctions of the possible shape optimizers for the third and fourth INE.
\begin{figure}
\centering
\begin{subfigure}{}
\includegraphics[width=5.5cm]{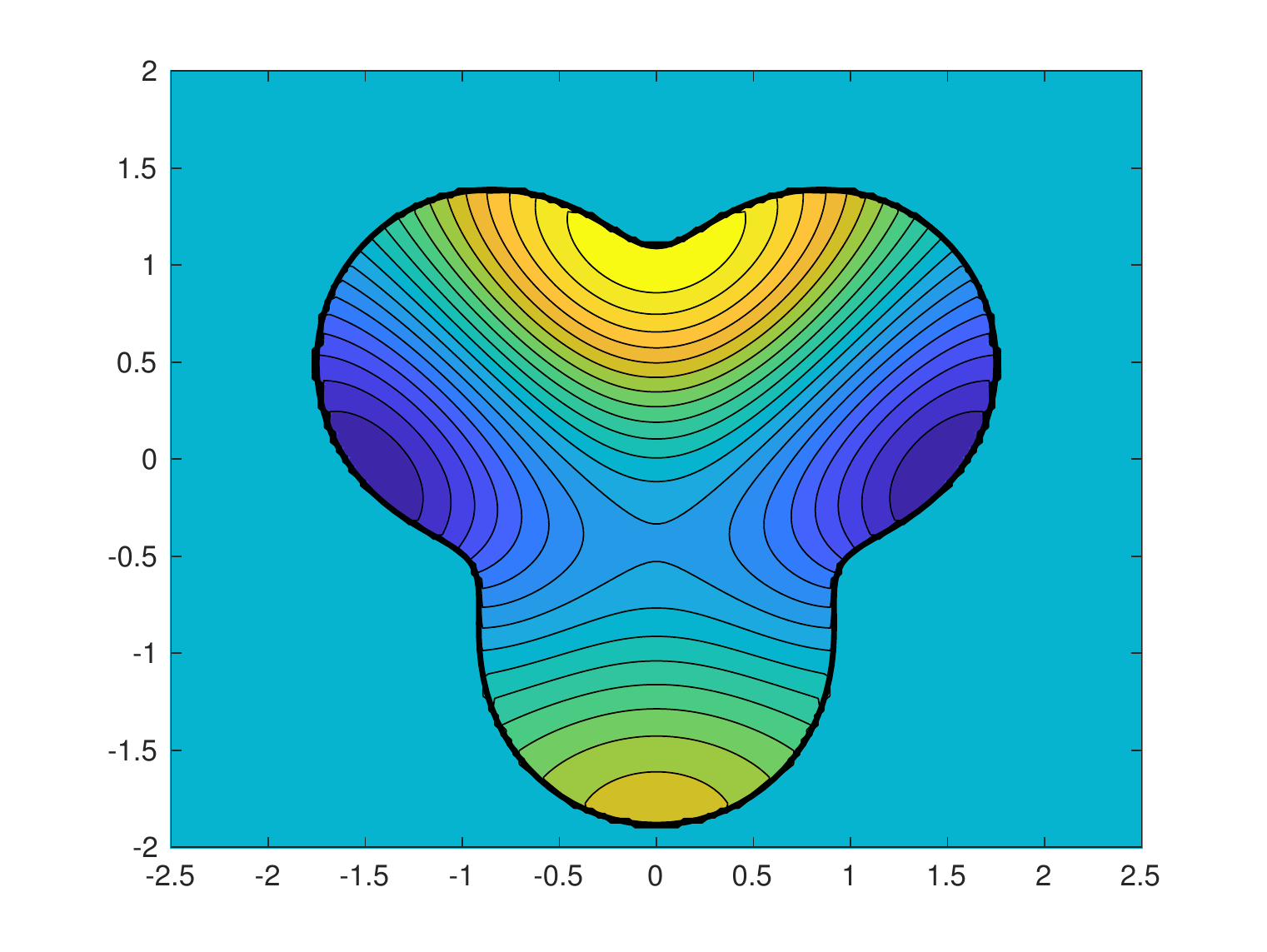}
\end{subfigure}\quad
\begin{subfigure}{}
\includegraphics[width=5.5cm]{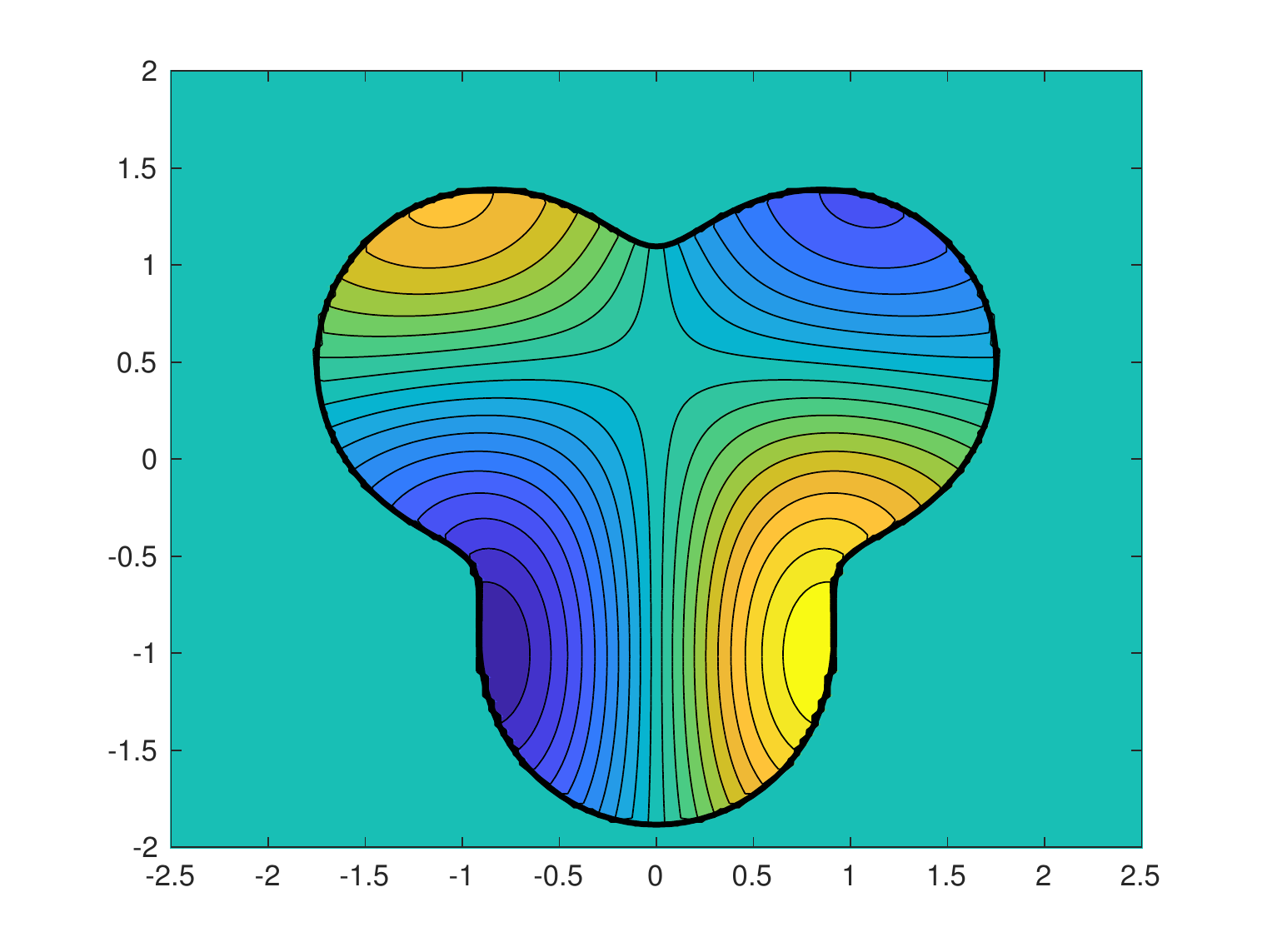}
\end{subfigure}\quad
\begin{subfigure}{}
\includegraphics[width=5.5cm]{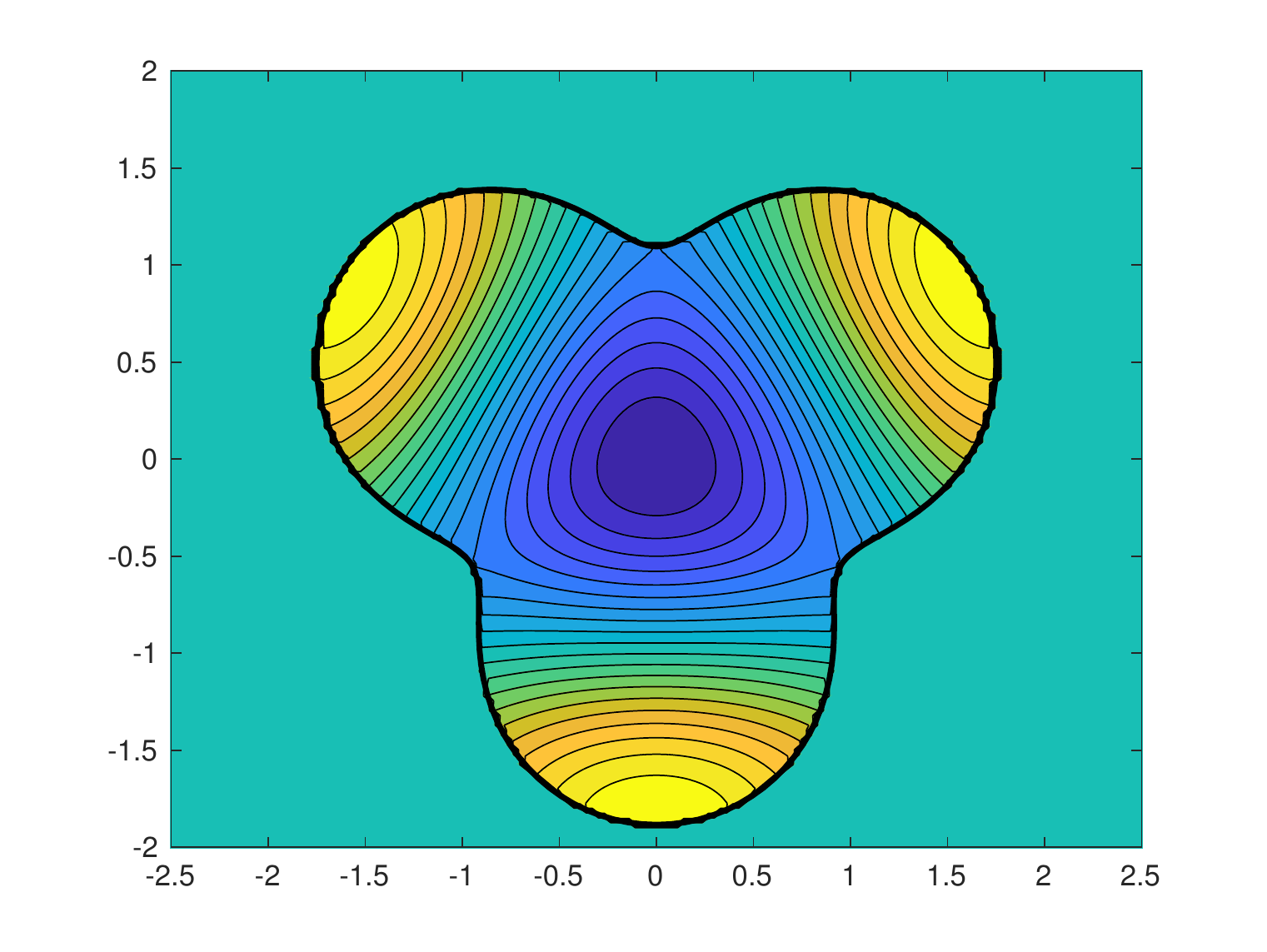}
\end{subfigure}\quad
\begin{subfigure}{}
\includegraphics[width=5.5cm]{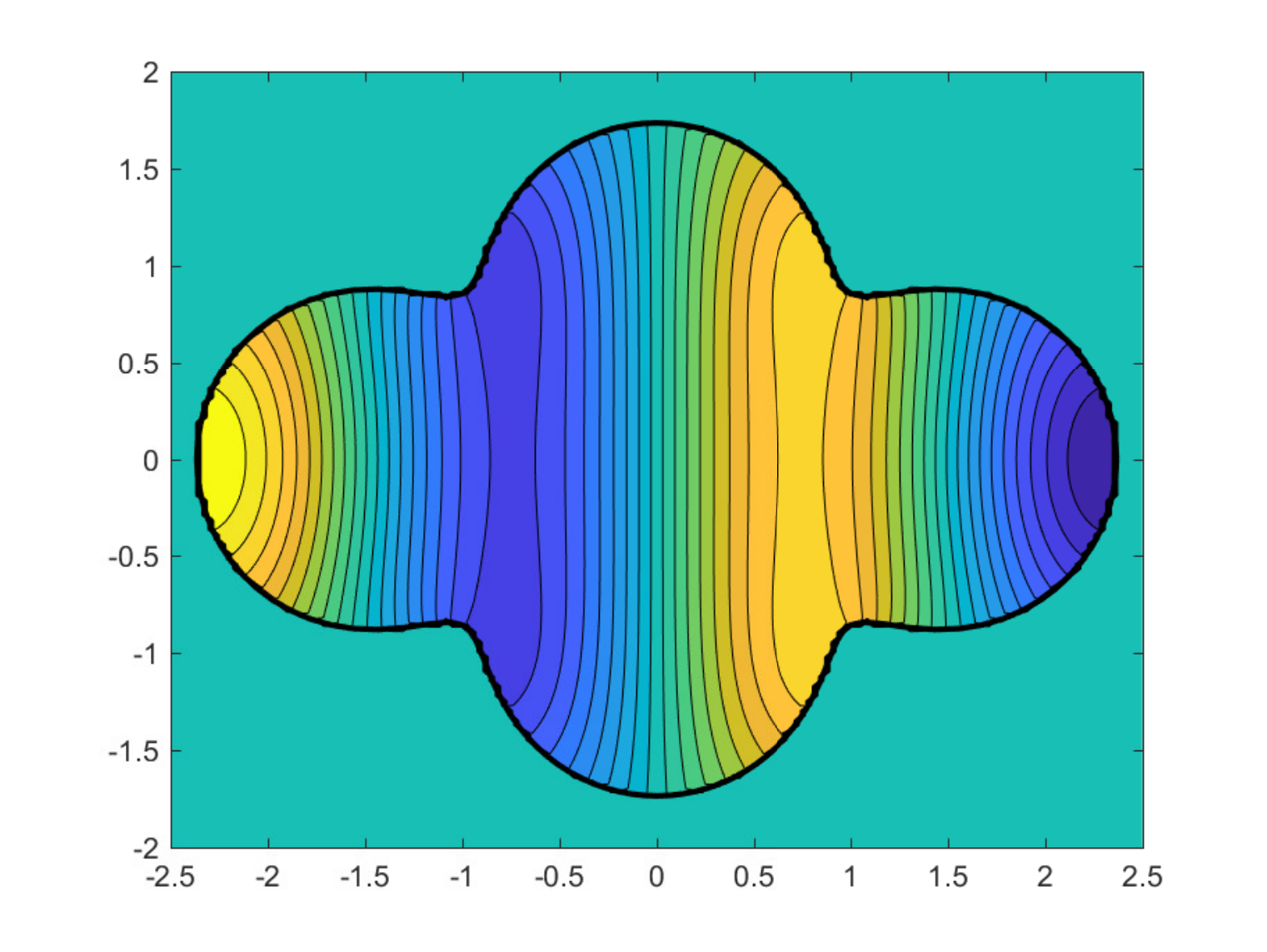}
\end{subfigure}\quad
\begin{subfigure}{}
\includegraphics[width=5.5cm]{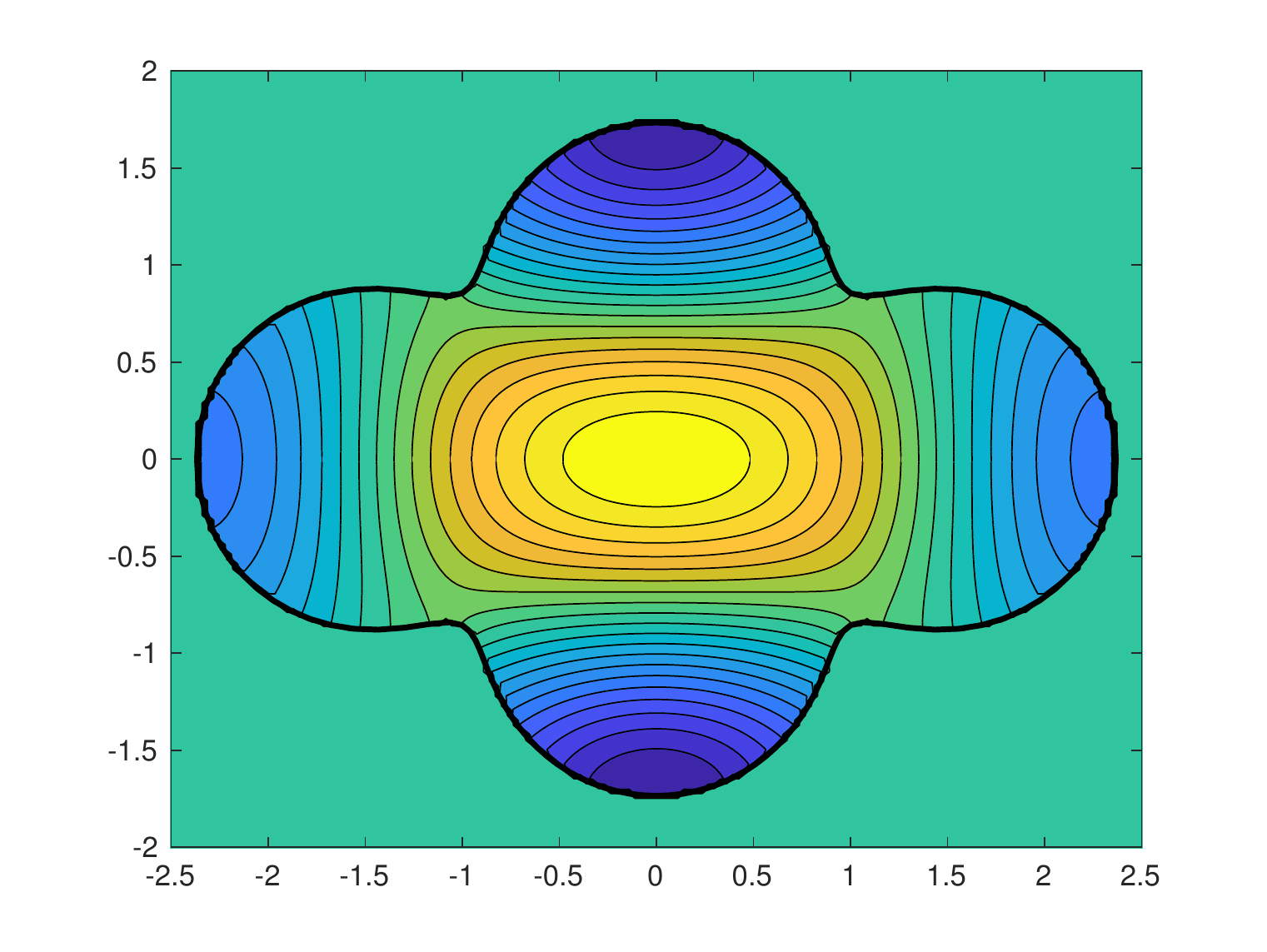}
\end{subfigure}\quad
\begin{subfigure}{}
\includegraphics[width=5.5cm]{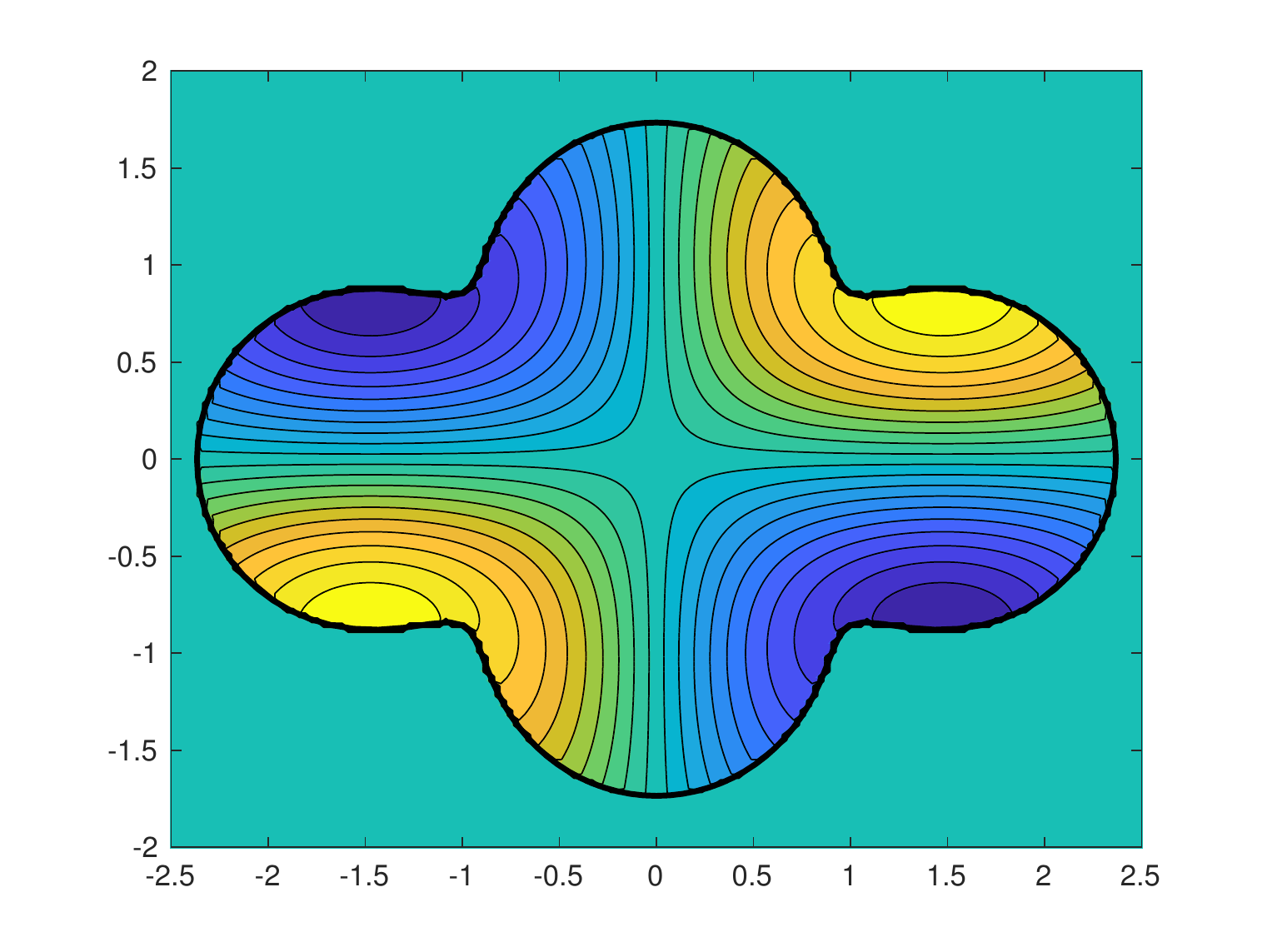}
\end{subfigure}
\caption{The three eigenfunction of the shape optimizer for the third and fourth INE. The parameters are $\alpha=2.0171$ and $c=1.6883$ with $32.9018$ having multiplicity three for the third INE and $\alpha=2.5426$ and $c=2.0845$ with $43.8694$ having multiplicity three.}
\label{fig:kle5}
\end{figure}

Note that we used $n=512$ for all numerical calculation to ensure that we have at least six digits accuracy for the values $\lambda_k\cdotp A$. This is guaranteed since we almost have a convergence of order four due 
to the fact that we have approximated the boundary and the unknown density function by quadratic interpolation (refer to \cite{KlLi12} for a superconvergence proof for three-dimensional scattering objects).

\section{Shape optimization for interior transmission eigenvalues}\label{sec:kle3}
Recall that interior transmission eigenvalues (ITEs) are numbers $\lambda=\kappa^2\in \mathbb{C}\backslash \{0\}$ such that 
    \begin{align*}
    \Delta w+\kappa^2 n w&=0\quad\;\;\text{in}\; \Omega\;,\\
    \Delta v+\kappa^2\;\;\; v&=0\quad\;\;\text{in}\;\Omega\;,\\
                           v&=w\quad\;\text{on}\;\partial\Omega\;,\\
            \partial_{\nu}v&=\partial_{\nu}w\;\text{on}\;\partial\Omega\;,
    \end{align*}
has a non-trivial solution $(v,w)\neq(0,0)$. Here, $n$ is the given index of refraction. This is a non-elliptic and non-self-adjoint problem. 
Existence and discreteness for real-valued $\kappa$ has already been established. However, the existence is still 
open for complex-valued $\kappa$ except for special geometries. To compute such ITEs for a given shape 
is therefore very challenging. 
We use the same technique as presented before for the numerical calculation of interior Neumann eigenvalues; that is, reduce the problem to a 
system of boundary integral equations, discretize it via a boundary element collocation method, and solve the resulting non-linear eigenvalue problem via the method of Beyn (see \cite{Be12}). 
For more details, we refer the reader to \cite{Kl13, Kl15} where ITEs for three-dimensional domains 
are computed and to \cite{KlPi18} for a good introduction for other methods to compute such ITEs.
Straightforwardly looking at real-valued ITEs using the index of refraction $n=4$ for different domains taken from \cite{KlPi18} reveals that neither the circle is maximizing nor minimizing $\lambda_1=A\cdotp \kappa_1^2$. 
The values $\lambda_1$ for eight different domains are given in Fig. \ref{fig:kle6}
\begin{figure}
\centering
\begin{subfigure}{}
\includegraphics[width=2.5cm]{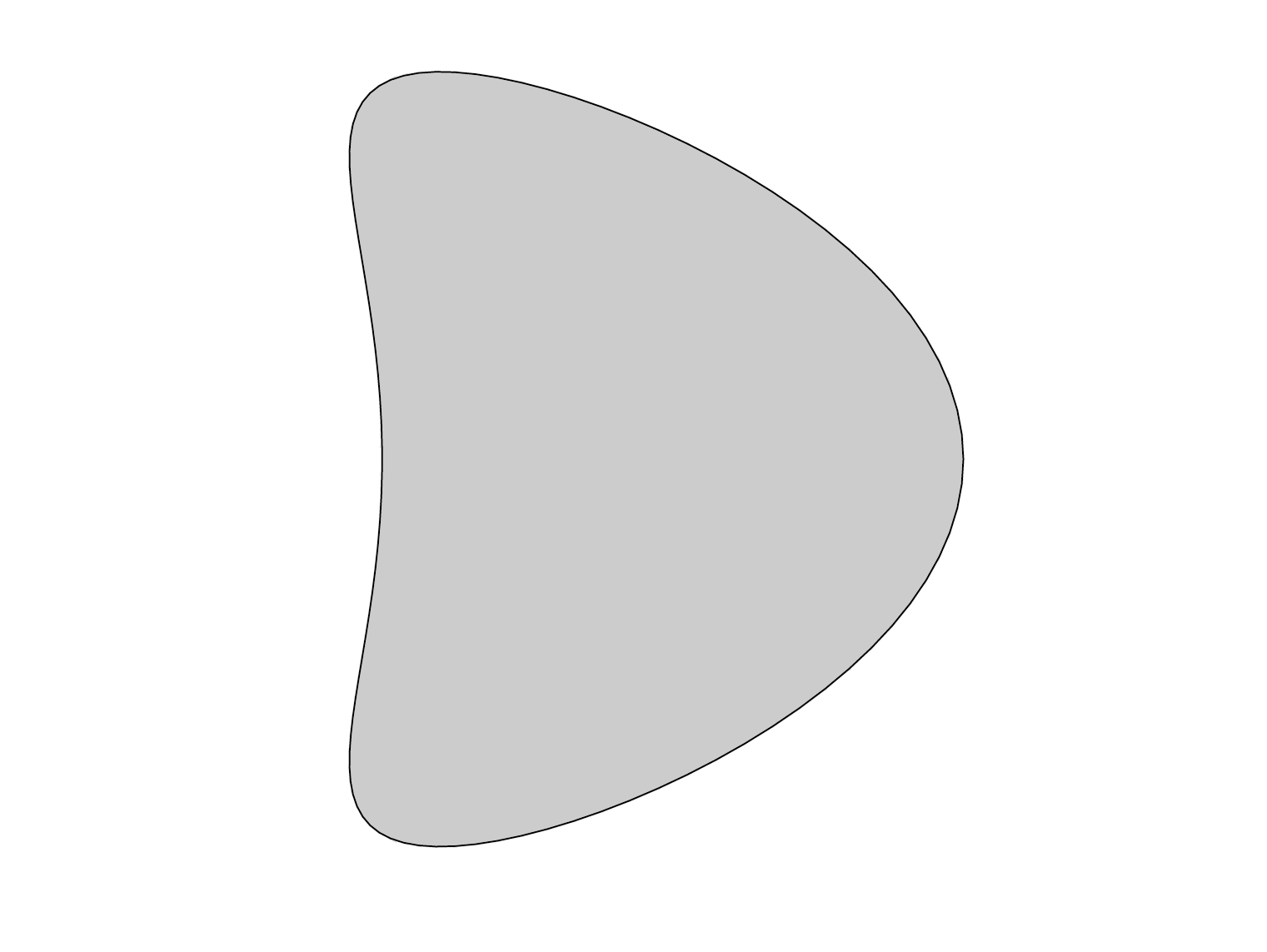}
\end{subfigure}\quad
\begin{subfigure}{}
\includegraphics[width=2.5cm]{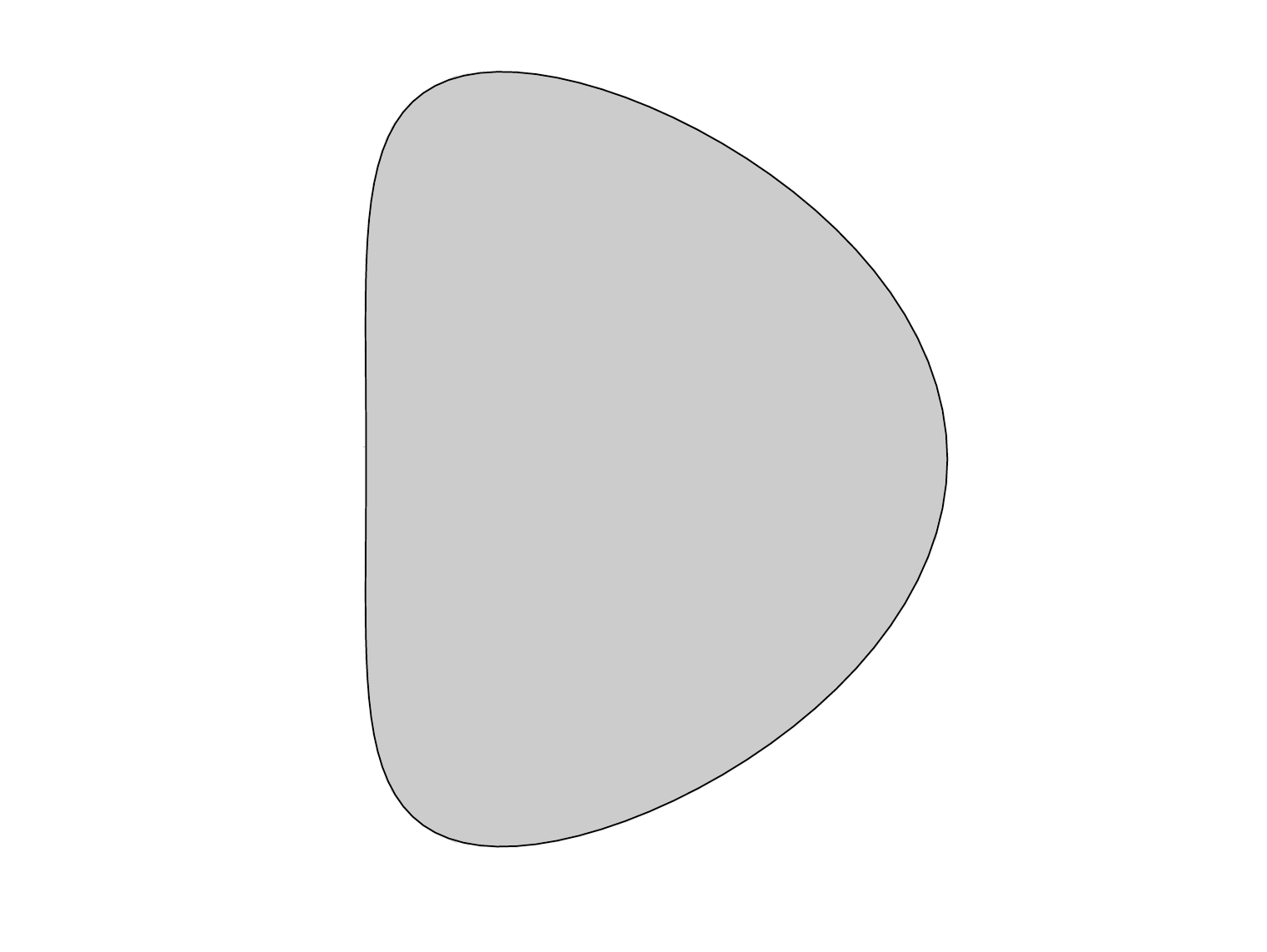}
\end{subfigure}\quad
\begin{subfigure}{}
\includegraphics[width=2.5cm]{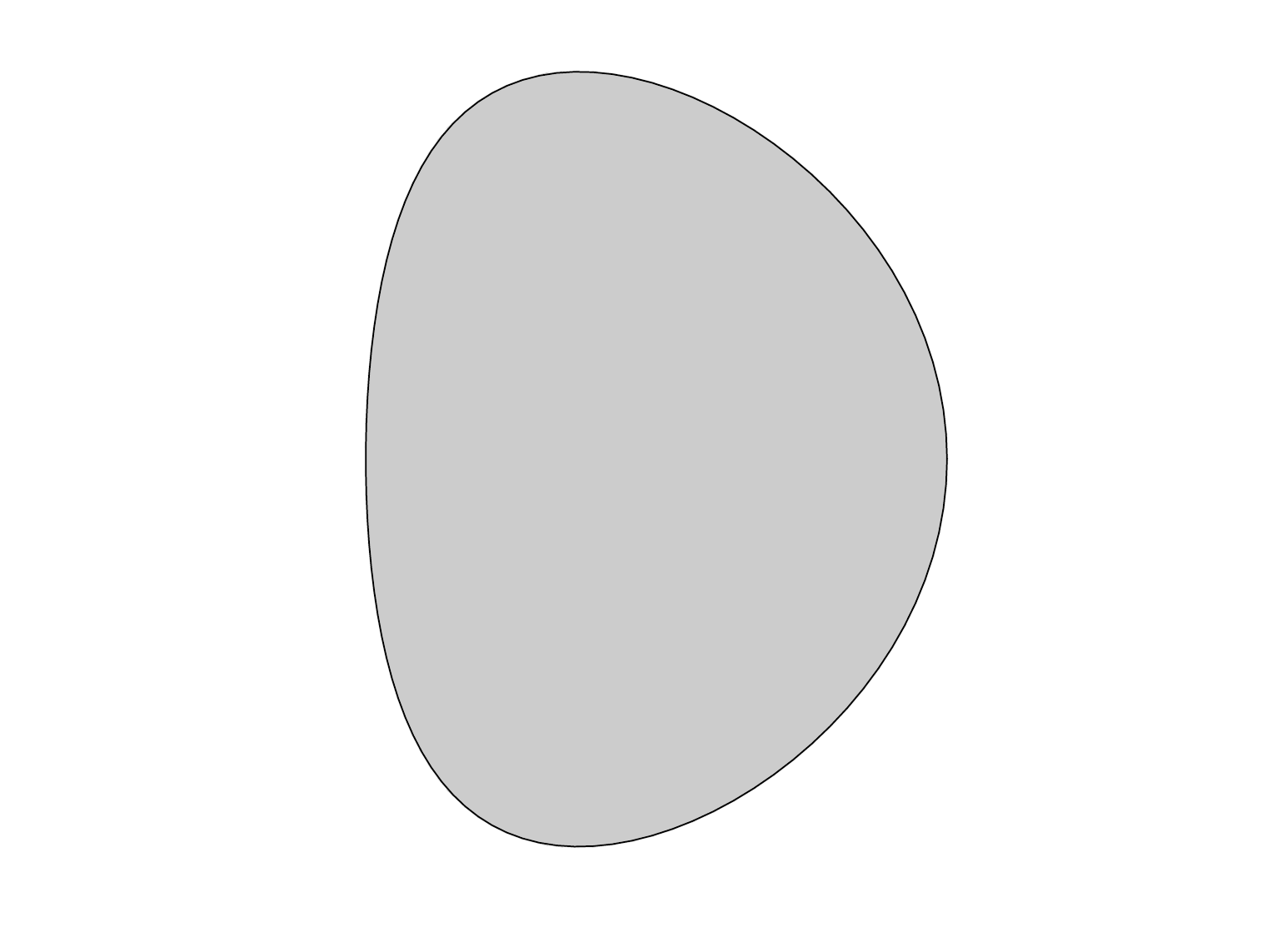}
\end{subfigure}\quad
\begin{subfigure}{}
\includegraphics[width=2.5cm]{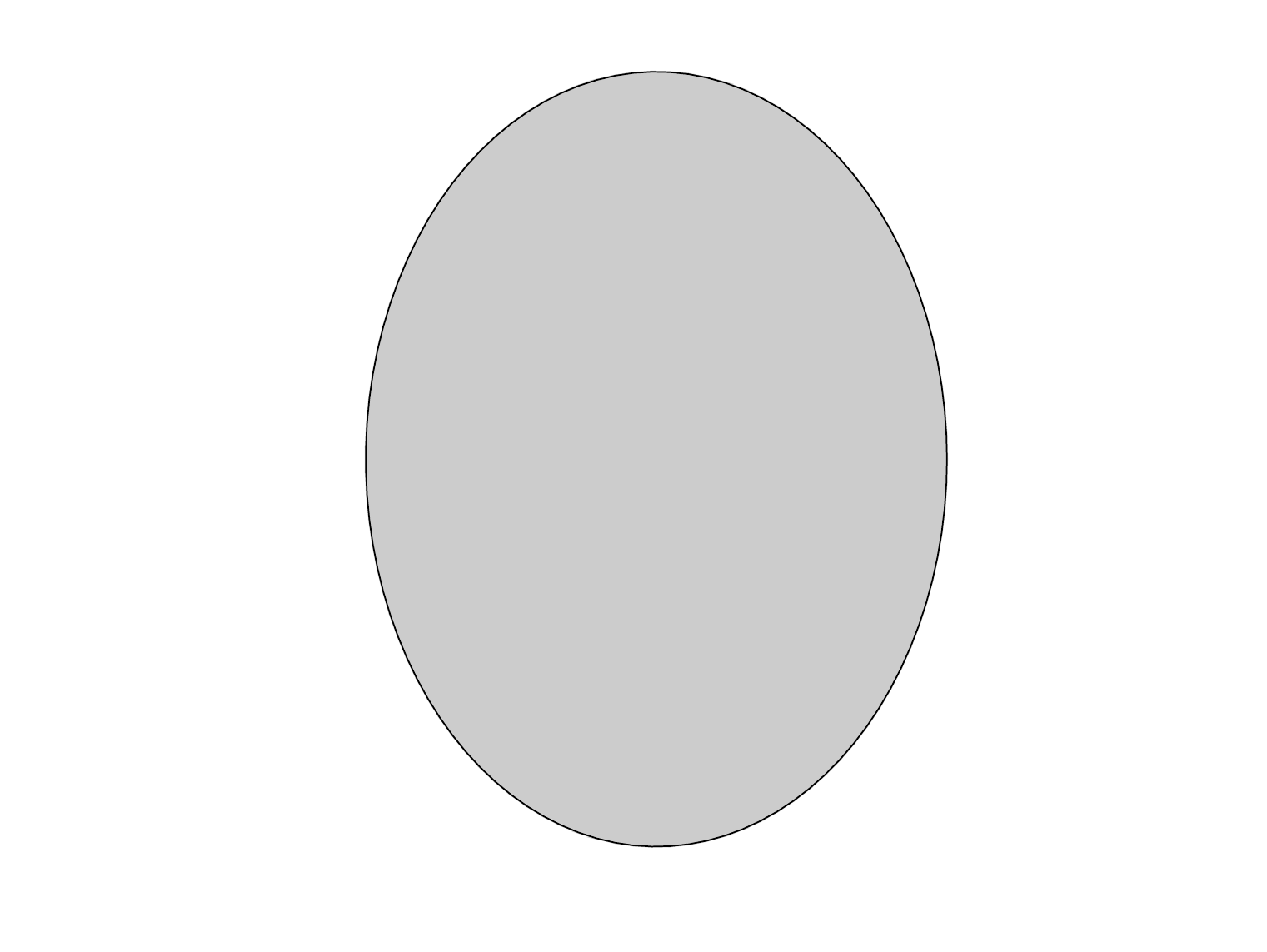}
\end{subfigure}\\
$29.1348\qquad\qquad\qquad\quad 26.9563\qquad\qquad\qquad\quad 25.2925\qquad\qquad\qquad\quad 24.6688$\\[6pt]
\begin{subfigure}{}
\includegraphics[width=2.5cm]{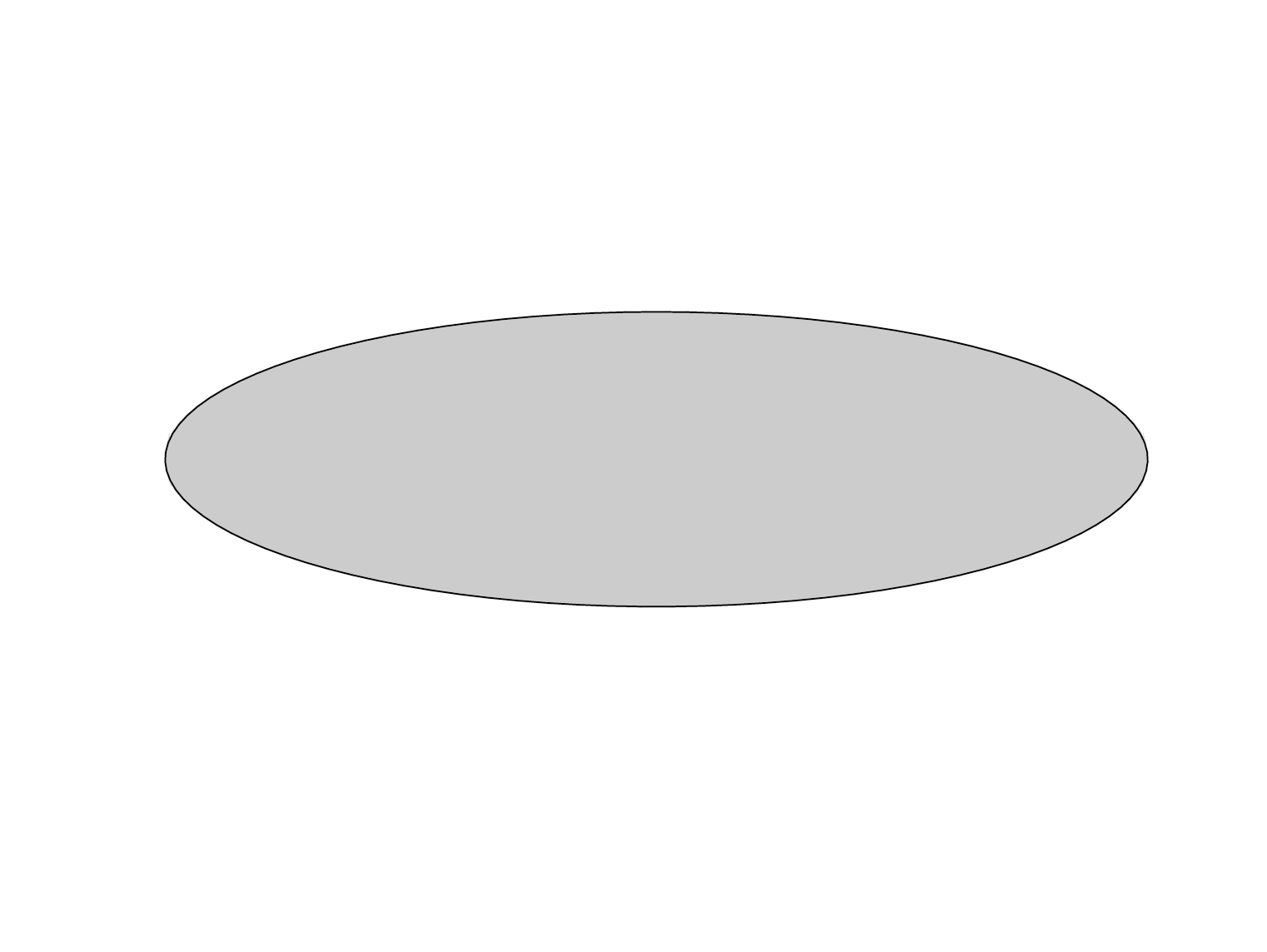}
\end{subfigure}\quad
\begin{subfigure}{}
\includegraphics[width=2.5cm]{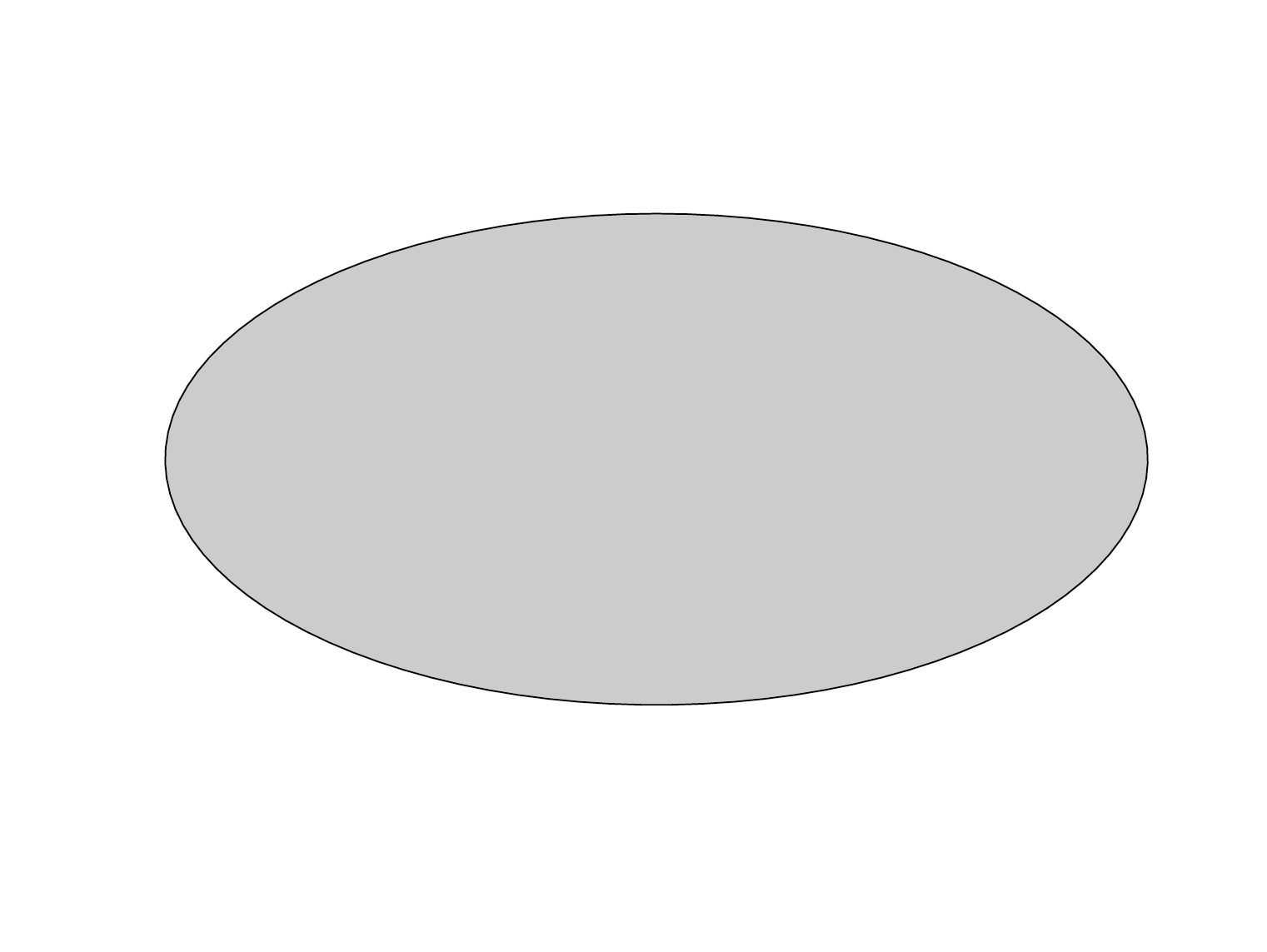}
\end{subfigure}\quad
\begin{subfigure}{}
\includegraphics[width=2.5cm]{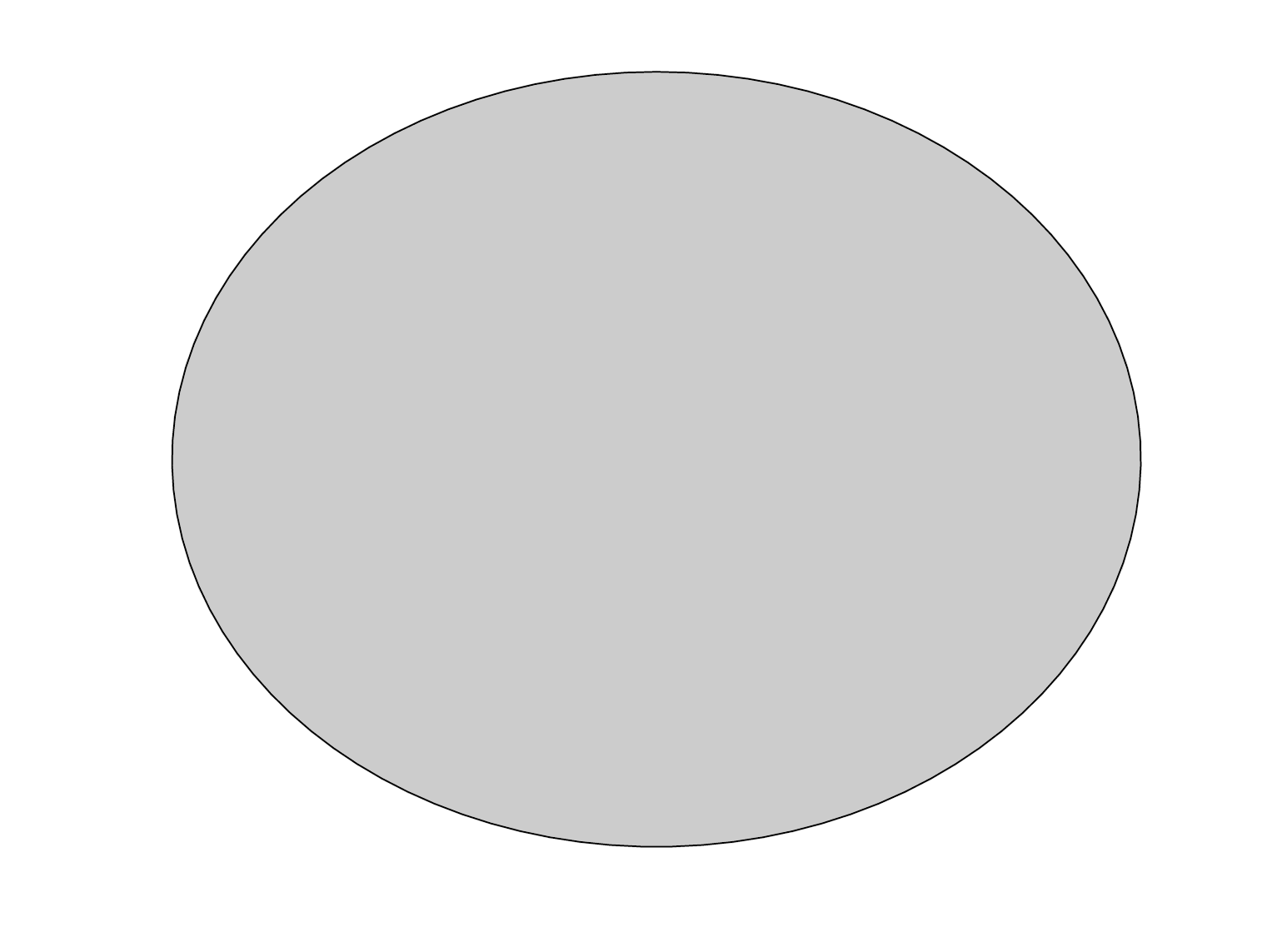}
\end{subfigure}\quad
\begin{subfigure}{}
\includegraphics[width=2.5cm]{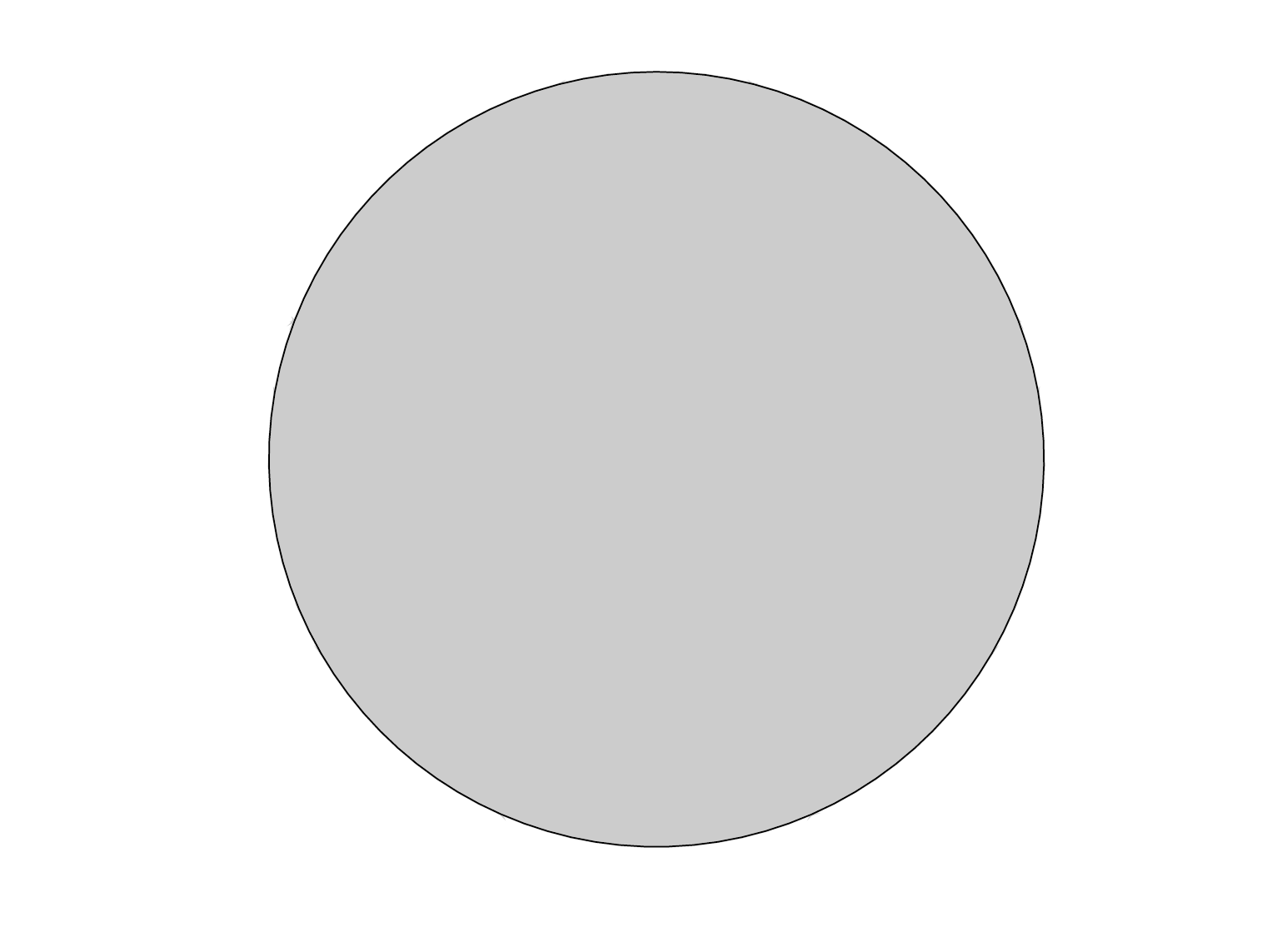}
\end{subfigure}\\
$40.4687\qquad\qquad\qquad\quad 29.4600\qquad\qquad\qquad\quad 24.7064\qquad\qquad\qquad\quad 26.4683$\\[6pt]
\caption{The values $\lambda_1$ for eight different domains using $n=4$.}
\label{fig:kle6}
\end{figure}

But recall that there might be complex-valued ITEs as well which are not taken into account. If we consider $|\lambda_1|$ instead of $\lambda_1$ using the same eight domains,
we obtain the results as presented in Fig. \ref{fig:kle7}.
\begin{figure}
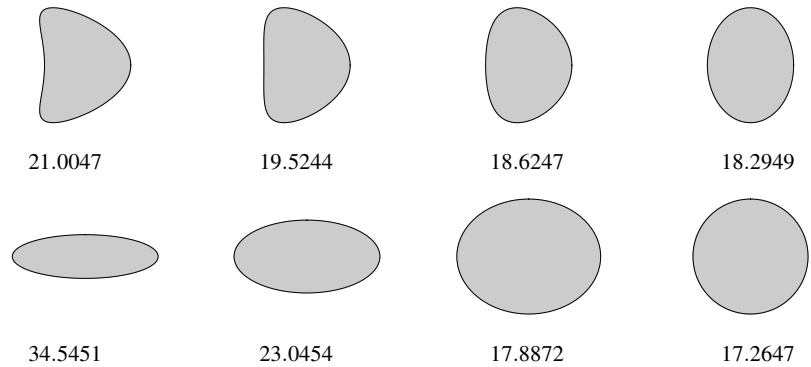

\centering
\begin{subfigure}{}
\includegraphics[width=2.5cm]{kleefeld-defo03N-eps-converted-to.pdf}
\end{subfigure}\quad
\begin{subfigure}{}
\includegraphics[width=2.5cm]{kleefeld-defo02N-eps-converted-to.pdf}
\end{subfigure}\quad
\begin{subfigure}{}
\includegraphics[width=2.5cm]{kleefeld-defo01N-eps-converted-to.pdf}
\end{subfigure}\quad
\begin{subfigure}{}
\includegraphics[width=2.5cm]{kleefeld-defo00N-eps-converted-to.pdf}
\end{subfigure}\\
$21.0047\qquad\qquad\qquad\quad 19.5244\qquad\qquad\qquad\quad 18.6247\qquad\qquad\qquad\quad 18.2949$\\[6pt]
\begin{subfigure}{}
\includegraphics[width=2.5cm]{kleefeld-elli03N-eps-converted-to.pdf}
\end{subfigure}\quad
\begin{subfigure}{}
\includegraphics[width=2.5cm]{kleefeld-elli05N-eps-converted-to.pdf}
\end{subfigure}\quad
\begin{subfigure}{}
\includegraphics[width=2.5cm]{kleefeld-elli08N-eps-converted-to.pdf}
\end{subfigure}\quad
\begin{subfigure}{}
\includegraphics[width=2.5cm]{kleefeld-elli10N-eps-converted-to.pdf}
\end{subfigure}\\
$34.5451\qquad\qquad\qquad\quad 23.0454\qquad\qquad\qquad\quad 17.8872\qquad\qquad\qquad\quad 17.2647$\\[6pt]
\caption{The values $|\lambda_1|$ for eight different domains using $n=4$.}
\label{fig:kle7}
\end{figure}

As one can observe, it seems that the circle is minimizing $|\lambda_1|$. Hence, if we consider 
$|\lambda_1|\leq |\lambda_2| \leq |\lambda_3|\leq \cdots$, then we make the \textbf{conjecture} that the first absolute ITE is minimal for a circle for the index of refraction $n>1$.
If this is true, then it is also true for $0<n<1$ using the relation $\kappa(1/n)=\sqrt{n}\kappa(n)$. Further, since $\lambda_1$ is complex-valued, it comes in complex conjugate pairs. 
Hence, the second eigenvalue will be minimized by a circle as well.

Further investigation of shapes that minimize higher interior transmission eigenvalues is a very interesting and challenging topic.

\section{Summary and outlook}\label{sec:kle4}
In this paper, it is shown how to efficiently compute interior Neumann eigenvalues for a given domain in two dimensions. 
Additionally, the value of the shape maximizer for the third and fourth interior Neumann eigenvalue has been improved from $32.90$ and $43.86$ to $32.9018$ and $43.8694$ with multiplicity three, respectively. 
At the same time, the number of parameters describing the boundary of a possible maximizer has been reduced to two parameters using modified equipotentials. 
The conjecture is that the third and fourth interior Neumann eigenvalue might be given by such modified equipotentials. 
This work presents very recent numerical results and a further investigation has to be carried out in order to validate whether the shape maximizer for higher interior Neumann eigenvalues can be found with modified equipotentials. This idea can easily be used for extending this approach to the three-dimensional case.

Moreover, for the first time numerical results are presented for the minimization of interior transmission eigenvalues in two dimensions although already the numerical calculation of those for a given domain is a 
very challenging task since the problem is neither elliptic nor self-adjoint and hence complex-valued interior transmission eigenvalues might exist. From the theoretical point of view, this fact is still open. 
Additionally, it is open whether there exist a unique minimizer for the first and second interior transmission eigenvalue. Here, we show numerically and hence conjecture that the first and second interior 
transmission eigenvalue is minimized by a circle. It remains to prove this observation, but it cannot be carried out by standard spectral arguments like for the Dirichlet, Neumann, Robin, or Steklov 
eigenvalue problem. Moreover, one can now try to investigate the three-dimenensional case.

Above all, one could also investigate the electromagnetic and/or the elastic scattering case in two and three dimensions.
        
\section*{Acknowledgement}
I would like to thank the IMSE'18 steering committee for giving me the opportunity to present my recent results for the maximization of interior Neumann  
and minimization of interior transmission eigenvalues on July 19th, 2018. Further, I would like to thank Paul Harris for the organization of this nice event at the University of Brighton, UK.

\par\egroup

\printindex   
\end{document}